\documentclass[12pt,oneside,reqno]{amsart}

\title{Quasi-Hamiltonian Model Spaces}
\author{Kay Paulus}
\address{K. Paulus}
\email{math@kaypaulus.de}

\author{Bart Van Steirteghem}
\address{B. Van Steirteghem:  Friedrich-Alexander-Universit\"at Erlangen-N\"urnberg}
\email{bartvs@math.fau.de}

\usepackage[english]{babel}

\usepackage{fullpage}
\usepackage{amssymb, amsmath}
\usepackage{mathrsfs}
\usepackage[shortlabels]{enumitem}
\usepackage[colorlinks,breaklinks, allcolors=blue]{hyperref}
\usepackage{longtable}
\usepackage{cleveref}

\usepackage{dynkin-diagrams}
\pgfkeys{/Dynkin diagram, 
edge length=.5cm,
}

\swapnumbers
\newtheorem{thm}{Theorem}[section]
\newtheorem{lem}[thm]{Lemma}
\newtheorem{prop}[thm]{Proposition}
\newtheorem{cor}[thm]{Corollary}

\theoremstyle{definition}
\newtheorem{defn}[thm]{Definition}

\newtheorem{rem}[thm]{Remark}

\numberwithin{equation}{section}

\newcommand{\R}{\mathbb{R}}
\newcommand{\Z}{\mathbb{Z}}
\newcommand{\N}{\mathbb{N}}
\newcommand{\K}{\mathbb{K}}
\newcommand{\Q}{\mathbb{Q}}
\newcommand{\C}{\mathbb{C}}

\newcommand{\Chi}{\mathcal X}

\newcommand{\SL}{\mathrm{SL}}
\newcommand{\Spin}{\mathrm{Spin}}

\newcommand{\Sp}{\mathrm{Sp}}
\newcommand{\SU}{\mathrm{SU}}

\newcommand{\fa}{\mathfrak{a}}

\newcommand{\fk}{\mathfrak{k}}

\newcommand{\ft}{\mathfrak{t}}

\newcommand{\sv}{\mathsf{v}}

\DeclareMathOperator{\Hom}{Hom}
\DeclareMathOperator{\Lie}{Lie}

\DeclareMathOperator{\Ad}{Ad}

\newcommand{\lat}{\Lambda}
\newcommand{\wm}{\Gamma}
\newcommand{\dw}{\wl_+}
\newcommand{\wl}{P} 
\newcommand{\w}{\omega}
\newcommand{\wt}{\widetilde{\omega}}

\newcommand{\Kt}{K\tau}
\newcommand{\Kvt}{K_{\mathsf{v}\tau}}
\newcommand{\Kvtc}{K_{\mathsf{v}\tau}^{\C}}
\newcommand{\ww}{\tilde\w}
\newcommand{\Rp}{\R_{\ge 0}}
\newcommand{\ofa}{\overline{\fa}}
\newcommand{\oS}{\overline{S}}
\newcommand{\rank}{\operatorname{rank}}
\newcommand{\la}{\langle}
\newcommand{\ra}{\rangle}
\newcommand{\SN}{\Sigma^N}

\renewcommand{\a}{\alpha}
\renewcommand{\aa}{\overline{\alpha}}
\newcommand{\aav}{\aa^{\vee}}
\newcommand{\bet}{\overline{\beta}}
\renewcommand{\P}{\mathcal{P}}
\newcommand{\A}{\mathcal{A}}
\newcommand{\G}{\Gamma}

\begin{document}

\begin{abstract}
Let $K$ be a simple and simply connected compact Lie group. We call a  (twisted) quasi-Hamiltonian $K$-manifold $M$ a quasi-Hamiltonian model space if it is multiplicity free and its momentum map is surjective. We explicitly identify the subgroups of the Lie algebra of the maximal torus of $K$, which, by F.~Knop's  classification of multiplicity free quasi-Hamiltonian manifolds, are in one-to-one correspondence with the isomorphism classes of quasi-Hamiltonian model $K$-spaces. 
\end{abstract}

\maketitle

\section{Introduction}
A quasi-affine variety equipped with an action of a complex connected reductive group $G$ is called a \emph{model variety} for $G$ if its coordinate ring contains every irreducible representation of $G$ exactly once. The study of such `representation models' started in \cite{bgg-model-russian} and has been quite fruitful, see for example \cite{gelfand-zelevinski-model1, gelfand-zelevinski-model2, adams-huang-vogan, luna-model,bravi_gandini_maffei-proj_norm_model}.

In this paper, we classify analogous \emph{model spaces} in the setting of the quasi-Hamiltonian manifolds introduced by A.~Alekseev, A.~Malkin and E.~Meinrenken in \cite{alekseev-malkin-meinrenken-groupvaluedmoment}. 
Roughly speaking, a quasi-Hamiltonian $K$-manifold is a smooth manifold $M$ equipped with an action of a compact connected Lie group $K$, a 2-form $\w$ and a smooth $K$-equivariant map $m:M\to K$, called the {\em (group valued) momentum map}, fulfilling certain compatibility conditions (see \cref{def:qHam}).

In fact, this notion can be generalized by allowing a twist of the conjugation action of $K$ on itself. Indeed, given a smooth automorphism $\tau$ of $K$ one can require that the momentum map $m: M \to K$ be equivariant with respect to the twisted conjugation action \[k\cdot_\tau g = kg\tau(k)^{-1}\] of $K$ on itself.  In this case, we use $K\tau$ for $K$ equipped with this $\tau$-twisted action and denote the momentum map by $m:M \to K\tau$. Such quasi-Hamiltonian $\Kt$-manifolds were first defined by Meinrenken in \cite{meinrenken-twistedconvex} and  independently by F.~Knop in \cite{knop-qham-arxiv}. In \cite{boalch_and_yamakawa}, P.~Boalch and D.~Yamakawa also independently considered such manifolds in the context of twisted wild character varieties.

From now on, \emph{we assume that $K$ is simply connected}. As is known and will be recalled in \cref{thm:twistconj}, there is a natural homeomorphism $c: \A \to K\tau/K$ between the set $K\tau/K$ of $\tau$-twisted conjugacy classes in $K$ and the fundamental alcove $\A$ of a certain affine root system.  In \cite[Theorem 7.2]{alekseev-malkin-meinrenken-groupvaluedmoment} and   \cite[Theorem 4.4]{meinrenken-twistedconvex} it was shown (for $\tau = id_K$ and general $\tau$, respectively) that when $M$ is a compact and connected quasi-Hamiltonian $\Kt$-manifold, then the image of the map 
\begin{equation}
m_+:=c^{-1} \circ \pi \circ m:M \to \A, \label{eq:mplus}
\end{equation}
where $\pi: K\tau \to K\tau/K$ is the quotient map, is a convex polytope $\P_M$ in $\A$, which is called the \textbf{momentum polytope} of $M$. 

Alekseev, Malkin and Meinrenken also ported the classical notion of symplectic reduction of Hamiltonian manifolds to the setting of quasi-Hamiltonian manifolds \cite[Section 5]{alekseev-malkin-meinrenken-groupvaluedmoment}. In analogy to the multiplicity free Hamiltonian manifolds of Guillemin and Sternberg \cite{guill&stern-mf}, Knop then made the following definition in \cite{knop-qham-arxiv}: a compact connected quasi-Hamiltonian manifold is called  \textbf{multiplicity free} if all its symplectic reductions are points, see also \cite[Def.\ 2.4.1 and Prop.\ 2.4.2]{knop-cmfqham-arxiv}.

In \cite[Corollary 2.6.2]{knop-cmfqham-arxiv}, Knop showed that compact connected multiplicity free quasi-Hamiltonian $\Kt$-manifolds $M$ are uniquely determined by the pair $(\P_M,\lat_M)$, where $\lat_M$ is a certain lattice  which encodes the principal isotropy group of the $K$-action on $M$. In addition, he characterized which pairs $(\P,\lat)$ consisting of a polytope and a lattice can occur this way.

Knop also studied certain series of examples of multiplicity free quasi-Hamiltonian manifolds in \cite[Section 2.7]{knop-cmfqham-arxiv}, and in \cite{paulus_phd} the first author obtained a classification of those for which $\dim \P_M=1$, see also \cite{knop-paulus-arxiv}.  In \cite[Proposition 2.7.2]{knop-cmfqham-arxiv}, Knop identified some multiplicity free quasi-Hamiltonian manifolds which are `as big as possible.' Their explicit combinatorial classification is the purpose of this paper.

\begin{defn} \label{def:modelspace}
A compact connected quasi-Hamiltonian $K\tau$-manifold is called a {\bf (quasi-Hamiltonian) model $K\tau$-space} if it is multiplicity free and its momentum map is surjective.
\end{defn}

The main result of this paper is \cref{thm_main} in which we combinatorially classify model $\Kt$-spaces for $K$ \emph{simple} and simply connected. The necessary prerequisites for stating this theorem are reviewed in \cref{arsc}. It is an application of \cref{sphpair}, which is a specialization of Knop's aforementioned classification theorem.  Since the momentum polytope of a model $\Kt$- space is always the full alcove $\A$, our  classification is in terms of the possible lattices $\lat_M$. It would be interesting to have global descriptions of the model spaces $M$ realizing these lattices.

Knop's characterization in \cite{knop-cmfqham-arxiv} of the pairs $(\P_M, \lat_M)$ realized by multiplicity free quasi-Hamiltonian manifolds $M$ is in terms of \emph{weight monoids} of smooth affine \emph{spherical varieties}. This weight monoid is a basic representation-theoretic invariant of such varieties (see \cref{def:sasv}).   In \cref{secloc} we present the  tools from the combinatorial theory of spherical varieties that we will use  in \cref{seclat} to prove \cref{thm_main}. The main tool is \cref{pvsalg}. It is a special case, adapted to our setting, of the combinatorial characterization of the weight monoids of smooth affine spherical varieties in \cite{charwm} and may be of some independent interest.

\subsection*{Notation}
Unless stated otherwise, $K$ will be a simple  and simply connected compact Lie group with Lie algebra $\mathfrak{k}$. Furthermore $K$ will be equipped with a (possibly trivial) smooth automorphism $\tau$, which we will also call a twist, and its Lie algebra $\fk$ with a scalar product $\la\cdot,\cdot \ra$ which is invariant for $K$ and $\tau$. 
When $A$ is a subset of a free abelian group $\Chi$, we will use $\Z A$ for the smallest subgroup of $\Chi$ containing $A$ and when $A=\{a_1,a_2,\ldots,a_n\}$ is a finite set, we will also use  $\la a_1,a_2,\ldots,a_n\ra_{\Z}$ for this group.  When $D$ is a subset of a real vector space $V$, we will use $\Rp D$ for the closed convex cone generated by $D$ in $V$. 

\subsection*{Acknowledgment}
We thank Friedrich Knop for proposing the problem addressed in this paper, and for numerous helpful conversations. Part of this paper is based on the first author's doctoral thesis \cite{paulus_phd}, which was written under Knop's supervision. We also thank Guido Pezzini and Wolfgang Ruppert for many helpful discussions and Franziska Pechtl for her help with proofreading. Finally, we thank the referees of an earlier version of this paper for many helpful remarks and suggestions which led to improvements. The second author received support from the City University of New York PSC-CUNY Research Award Program.

\section{Prerequisites and main result}\label{arsc}
In this section we briefly recall, mostly following \cite{knop-cmfqham-arxiv}, the necessary notions to state both \cref{sphpair}, which is the special case of Knop's classification theorem \cite[Corollary 2.6]{knop-cmfqham-arxiv} that we will use, and \cref{thm_main}, which is our main result.

Although it will not play a direct role in what follows, we begin by giving, for completeness, the definition of a quasi-Hamiltonian $\Kt$-manifold, following \cite[Definition 2.1.2]{knop-cmfqham-arxiv}.

\begin{defn}\label{def:qHam}
  A {\bf quasi-Hamiltonian $K\tau$-manifold} is a smooth $K$-manifold
  $M$ equipped with a $K$-invariant $2$-form $\w$ and a
  $K$-equivariant smooth map $m:M\to K\tau$, called the (group valued)
  {\bf momentum map}, such that
 \begin{enumerate}
 \item  $d\w=-m^*\chi$,
 \item $\w(\xi x, \eta) = \la \xi, m^*\theta_{\tau}(\eta)\ra $ for all $\xi \in \fk, x \in M$ and $\eta \in T_xM$,
 \item $\ker \w_x = \{\xi x \in T_xM : \xi \in \fk \text{ with }\Ad m(x)({}^\tau \xi)+\xi = 0\}$,
 \end{enumerate}
 where $\theta_\tau:=\frac12({}^{\tau^{-1}}\theta^L+\theta^R)$ with
 $\theta^L, \theta^R$ being the left- and right-invariant
 Maurer-Cartan-forms on $K$ and 
\[\chi:=\frac{1}{12}\la \theta^L,[\theta^L, \theta^L]\ra = \frac{1}{12}\la \theta^R,[\theta^R, \theta^R]\ra \] 
is the canonical biinvariant closed $3$-form on $K$ with respect to the chosen scalar product $\la\cdot,\cdot\ra$ on $\fk$.
\end{defn}

We move on to  affine root systems, extracting from \cite[Section 1.1]{knop-cmfqham-arxiv}, which is based on \cite{macdonald-affine_rs_and_dedekind_eta} and \cite{macdonald-affine_Hecke_orthog_polys}, what we will need.  Let $\ofa$ be a Euclidean vector space  with inner product $(\cdot,\cdot)$ and  associated affine space $\fa$. We denote by $L(\fa)$ the set of affine linear functions on $\fa$. The gradient of $\a\in L(\fa)$ is denoted by $\aa \in \ofa$ and is characterized by the property
\begin{equation}\label{eq:afflin}
 \a(x+t)=\a(x)+(\aa, t),\ \text{ for all } x\in \fa, t\in \ofa.
\end{equation}

If $\a \in L(\fa)$ is a non-constant affine linear function,  
we define the \emph{reflection} 
\[s_{\a}: \fa \to \fa \text{ by } s_{\a}(x):= x-\a(x)\aa^{\vee}\]
where
\begin{equation} \label{eq:coroot}
\aa^{\vee} : = \frac{2}{(\aa,\aa)}\aa \in \ofa.
\end{equation}
Its induced action on an affine linear function $\beta\in L(\fa)$ is:
\[
 s_\a(\beta)=\beta-(\overline{\beta}, \aa^\vee) \a.
\]

\begin{defn}  \label{def_ars}
 A \textbf{(reduced) affine root system} on $\fa$ is a set $\Phi\subset L(\fa) \setminus\R$ of non-constant affine linear functions  such that:
 \begin{enumerate}[(a)]
  \item $\R\a\cap \Phi=\{\alpha, -\alpha\}$ for all $\a \in \Phi$,
  \item \label{def_ars_integral} $( \overline {\beta}, \aa^\vee) \in \Z$ for all $\a,\beta \in \Phi$
  \item $s_\a(\Phi)=\Phi$ for all $\a\in\Phi$,
  \item $\overline{\Phi}:=\{\aa \in \ofa : \a \in \Phi\}$ is finite.
 \end{enumerate}
 \end{defn}

Every $\a \in L(\fa) \setminus \R$ defines an affine hyperplane
\[H_{\a}:= \{x \in \fa: \a(x) = 0\}.\] An \textbf{alcove} of $\Phi$ is the closure of a connected component of $\fa \setminus \cup_{\a \in \Phi}H_\a$. The Weyl group $W_{\Phi}$ of $\Phi$, which is the subgroup generated by $\{s_\a : \a \in \Phi\}$ in the group of isometries of $\fa$, acts simply transitively on the set of alcoves of $\Phi$. Each such alcove $\A$ is a fundamental domain for the action of $W_\Phi$ on $\fa$.  

Put \[\overline\Phi:=\{\aa:\a\in\Phi\}\] and
$\overline{\Phi}^\vee:=\{\aa^\vee:\a\in\Phi\}$. These are possibly
non-reduced finite root systems on $\ofa$. 
 \begin{defn}An {\bf integral root system} on $\fa$ is a pair
   $(\Phi, \Xi)$ where $\Phi\subset L(\fa)$ is an affine root
   system and $\Xi\subseteq \ofa$ is a lattice with
   $\overline\Phi\subseteq \Xi$ and
   $(\Xi,\overline\Phi^\vee)\subseteq\Z$.
\end{defn}

Recall that $K$ is assumed to be simply connected. 
It is known that the twisted conjugacy classes in $K$ are in bijection with an alcove $\A$ of an affine root system that is determined by $K$ and  $\tau$, cf. \cite{wendt, mohrdieck-wendt-integralcc, meinrenken-twistedconvex}. We give the description of the twisted conjugacy classes from \cite[Section 2.2]{knop-cmfqham-arxiv}. To do so, we first recall that if $T \subset K$ is a maximal torus, then its character group $\Xi(T)$ can (and will) be identified with a sublattice of $\ft = \Lie(T)$ via the map
\[\Xi(T) \to \ft: \chi \mapsto a_{\chi},\]
where $a_{\chi}$ is the unique element of $\ft$ such that 
\[\chi(\exp \xi) = e^{2\pi i\la a_\chi,\xi\ra} \quad \text{for all } \xi \in \ft.\]
Consequently we also view the (finite) root system $\overline{\Phi}(\fk,\ft)$ of $K$ as a subset of $\ft$. In what follows, we will slightly abuse notation and no longer distinguish between $\fa$ and $\ofa$.  
\begin{thm}[{\cite[Theorem 2.2.1]{knop-cmfqham-arxiv}}] \label{thm:twistconj}
Let $\tau$ be an automorphism of the simply connected compact Lie group $K$. Then there exists a $\tau$-stable maximal torus $T$ in $K$ and an integral root system $(\Phi_{\tau}, \lat_{\tau})$ on the $\tau$-fixed part $\fa := \ft^{\tau}$ of $\ft$, with the following properties:
\begin{enumerate}[(a)]
\item If $\operatorname{pr}_\fa: \ft \to \fa$ is the orthogonal projection, then $\overline{\Phi}_\tau = \operatorname{pr}_\fa \overline{\Phi}(\fk,\ft)$ and $\lat_\tau = \operatorname{pr}_\fa \Xi(T)$. 
\item  \label{twistconjwl} The lattice $\lat_\tau$ is the weight lattice of $\overline{\Phi}_{\tau}$, that is  
\[\lat_\tau = \{\lambda \in \fa: \la \lambda, \aa^{\vee}\ra \in \Z \text{ for all } \a \in \Phi_{\tau}\}.\]
\item If $\A \subset \fa$ is an alcove of $\Phi_{\tau}$, then the composition
\(c: \A \subset \fa \stackrel{\exp}{\longrightarrow} K \rightarrow \Kt/K\)
is a homeomorphism. 
\item If for every $a \in \A$ we set
\begin{align*}
K_{a\tau}&:= \{ k \in K : k \cdot_{\tau} \exp(a) = \exp(a)\} \quad \text{and} \\
\Phi_\tau(a)&:=\{\a \in \Phi_\tau : \alpha(a)=0\}.
\end{align*}
then $K_{a\tau}$ is a closed connected subgroup of $K$ with maximal torus $\exp(a)$ and integral root system $\left(\overline{\Phi_\tau(a)},\lat_{\tau}\right)$.  
\end{enumerate}
Moreover, $(T, \Phi_{\tau})$ is unique up to twisted conjugation. 
\end{thm}

\begin{rem}
Using standard arguments (like those in the proof of \cite[Theorem 2.2.1]{knop-cmfqham-arxiv}, for example)  one shows that
\begin{enumerate}[(a)]
\item The type of the root system $\Phi_{\tau}$ in \cref{thm:twistconj} only depends on the image $\overline{\tau}$ of $\tau$ in the group of outer automorphisms of $K$;
\item If $K$ is simple then $\Phi_{\tau}$ is the irreducible affine root system of type $\mathsf{X}_n^{(r)}$ in \cref{ars}, where $\mathsf{X}_n$ is the Dynkin type of $K$ and  $r \in \{1,2,3\}$ is the order of $\overline{\tau}$. 
\end{enumerate}
\end{rem}

\begin{longtable}{c|c|c}
\caption{The Dynkin diagrams of the reduced and irreducible affine root systems, with Dynkin labels as given in \cite[Theorem 4.8]{kac-inf_dim_Lie_alg}. The Dynkin labels will play a role in \cref{seclat}.}\label{ars}\\
$\dynkin[labels*={1,1}]A[1]1$ &  \dynkin [labels={1},labels*={,1,1,1,1}, edge length=.8cm] A[1]{} &\begin{dynkinDiagram}[labels = {0,1,2,3,n-2,n-1,n}, label macro/.code={\alpha_{\drlap#1}},
 labels*={,,2,2,2,2,2},
  edge length=.8cm]B[1]{}
  \node[right,/Dynkin diagram/text style] at (root 0)
{\(1\)};
 \node[right,/Dynkin diagram/text style] at (root 1)
{\(1\)};
  \end{dynkinDiagram}
 \\ $\mathsf{A}_1^{(1)}$ & $\mathsf{A}_n^{(1)}, n\ge 2$ & $\mathsf{B}_n^{(1)}, n\ge 3$\\
 \hline
 \dynkin[labels={0,1,2,n-2,n-1,n}, label macro/.code={\alpha_{\drlap#1}}, labels*={1,2,2,2,2,1}, edge length=.7cm] C[1]{}&  
 \begin{dynkinDiagram}[labels={0,1,,3,n-3,,n-1,n}, label macro/.code={\alpha_{\drlap#1}}, labels*={,,2,2,2,2,,}, edge length=1cm] D[1]{}
 \node[ right,/Dynkin diagram/text style] at (root 0)
{\(1\)}; 
 \node[ right,/Dynkin diagram/text style] at (root 1)
{\(1\)}; 
 \node[ left,/Dynkin diagram/text style] at (root 2)
{\(\alpha_2\)}; 
 \node[ right,/Dynkin diagram/text style] at (root 5)
{\(\alpha_{n-2}\)}; 
 \node[ left,/Dynkin diagram/text style] at (root 6)
{\(1\)}; 
 \node[ left,/Dynkin diagram/text style] at (root 7)
{\(1\)}; 
 \end{dynkinDiagram}&\begin{dynkinDiagram}[labels*={,1,,2,,2,1}, label macro/.code={\alpha_{\drlap#1}}, label, edge length=.7cm] E[1]6 
  \node[above left,/Dynkin diagram/text style] at (root 4){\(3\)}; 
    \node[ left,/Dynkin diagram/text style] at (root 2){\(2\)}; 
      \node[ left,/Dynkin diagram/text style] at (root 0){\(1\)}; 
\end{dynkinDiagram} 
 \\
$\mathsf{C}_n^{(1)}, n \ge 2$ & $\mathsf{D}_n^{(1)}, n\ge 4$& $\mathsf{E}_6^{(1)}$\\
\hline
\begin{dynkinDiagram}[label, label macro/.code={\alpha_{\drlap#1}}, labels*={1,2,,3,4,3,2,1}, edge length=.7cm] E[1]7
 \node[left,/Dynkin diagram/text style] at (root 2){\(2\)};
\end{dynkinDiagram}
& \begin{dynkinDiagram}[label, label macro/.code={\alpha_{\drlap#1}}, labels*={1,2,,4,6,5,4,3,2}, edge length=.7cm] E[1]8
 \node[left,/Dynkin diagram/text style] at (root 2){\(3\)};
\end{dynkinDiagram}
&
 \dynkin[label, label macro/.code={\alpha_{\drlap#1}}, labels*={1,2,3,4,2},] F[1]4\\
$\mathsf{E}_{7}^{(1)}$ & $\mathsf{E}_{8}^{(1)}$ & $\mathsf{F}_4^{(1)}$\\
\hline
\dynkin[labels*={1,3,2}, ordering=Dynkin, label macro/.code={\alpha_{\drlap#1}}, label] G[1]2& \dynkin[label, label macro/.code={\alpha_{\drlap#1}}, labels*={2,1}] A[2]2 & \dynkin[labels = {0,1,2,3,n-2,n-1,n}, label macro/.code={\alpha_{\drlap#1}},
 label directions={,below}, labels*={2,2,2,2,2,2,1}, edge length=.7cm] A[2]{even} \\
 $\mathsf{G}_2^{(1)}$ & $\mathsf{A}_2^{(2)}$ & $\mathsf{A}_{2n}^{(2)}, n\ge 2$\\
 \hline
\begin{dynkinDiagram}[labels = {0,1,2,3,4,n-2,n-1,n}, label macro/.code={\alpha_{\drlap#1}}, labels*={1,1,,2,2,2,2,1}, edge length=.7cm, label* directions={,,above,,,,,}] A[2]{odd} 
 \node[above right,/Dynkin diagram/text style] at (root 2)
{\(2\)};
\end{dynkinDiagram}
& \dynkin[labels={0,1,2,n-2,n-1,n}, label macro/.code={\alpha_{\drlap#1}}, labels*={1,1,1,1,1,1}, edge length=.7cm] D[2]{}&
 \dynkin[label, label macro/.code={\alpha_{\drlap#1}}, labels*={1,2,3,2,1}, edge length=.5cm] E[2]6\\
$\mathsf{A}_{2n-1}^{(2)}, n\ge 3$ & $\mathsf{D}_{n+1}^{(2)}, n\ge 2$ & $\mathsf{E}_6^{(2)}$\\
\hline
&\dynkin[label, label macro/.code={\alpha_{\drlap#1}}, labels*={1,2,1}, edge length=.7cm] D[3]4 &\\
 &$\mathsf{D}_4^{(3)}$ &
\end{longtable}

Suppose now, that $M$ is a compact connected quasi-Hamiltonian $\Kt$-manifold, with momentum map $m$. Fixing an alcove $\A$ and a homeomorphism $c$ as in \cref{thm:twistconj}, one defines the so-called \emph{invariant momentum map} $m_+: M \to \A$ as in \cref{eq:mplus}. Recall from the introduction that its image
\[\P_M:= m_+(M) \subset \A\]
is a convex polytope, called the momentum polytope of $M$. 

The second invariant used in Knop's classifcation theorem of compact connected \emph{multiplicity free} quasi-Hamiltonian $\Kt$-manifolds is a lattice in $\fa$ which encodes the principal isotropy group of the $K$-action on $M$. We introduce it following \cite[Section 2.4]{knop-cmfqham-arxiv}. By \cref{thm:twistconj}, the isotropy goup $K_{a\tau}$  is the same subgroup of $K$ for all $a$ in the relative interior $\P^0_M$  of $\P_M$. Let's call this group $L_M$. Then the quotient of $L_M$ by  the kernel $L^1_M$ of its action on $m_+^{-1}(\P^0_M)$ is a torus which we call $A_M$. Furthermore $L^1_M$ is the principal isotropy group of the $K$-action on $M$ and it  is encoded by the character group $\lat_M$ of $A_M$, which we call the \textbf{lattice of $M$.} Because $\exp(\fa)$ is a maximal torus of $L_M$, the quotient map $L_M \to A_M$ restricts to a surjective homomorphism $\exp(a) \to A_M$. Consequently, the lattice $\lat_M$ can and will be viewed as subgroup of the lattice $\lat_{\tau}$, which is itself a subgroup of $\fa$.  

An immediate consequence of \cref{def:modelspace} is that the momentum polytope of a quasi-Hamiltonian model space $M$ is the alcove $\A$, so that only relevant invariant is the lattice $\lat_M$.
In order to characterize the lattices of quasi-Hamiltonian model spaces we need to make a few more recollections. 
Let $G$ be a complex connected reductive group and let $B \subset G$ be Borel subgroup. Write $\dw$ for the subset of $\Hom(B,\C^{\times})$ of dominant weights of $G$ (with respect to $B$). Recall that highest weight theory gives us a one-to-one correspondence $\lambda \to V(\lambda)$ between $\dw$ and the set of isomorphism classes of irreducible representations of $G$.  If $G$ acts on a variety $X$, then $G$ acts linearly on the ring $\C[X]$ of regular functions $X\to \C$ by
\[(g\cdot f)(x) := f(g^{-1}\cdot x) \quad \text{ for all $g\in G, f \in \C[X], x\in X$.}
\]
\begin{defn} \label{def:sasv}
A smooth affine irreducible $G$-variety is called \textbf{spherical} if its ring $\C[X]$ of regular functions is multiplicity free as a representation of $G$, that is, if every irreducible representation of $G$ occurs at most once in $\C[X]$. The \textbf{weight monoid} $\wm(X)$ of such a variety $X$ is the set of $B$-weights of $B$-eigenvectors in $\C[X]$, that is,
\[\wm(X) := \{\lambda \in \dw: \Hom_G(V(\lambda),\C[X]) \neq 0.\}.\]
\end{defn} 

\begin{rem}
Usually, a spherical $G$-variety is defined to be a normal $G$-variety $X$ which contains a dense orbit of the Borel subgroup $B$ of $G$. It follows from a well-known result due to Vinberg and Kimel'fel'd \cite{vin&kim} that the existence of a dense $B$-orbit on an affine $G$-variety $X$ is equivalent to the multiplicity-freeness of the $G$-module $\C[X]$. 
\end{rem}

Next we define the subgroups of $\fa$ that can occur as lattices of quasi-Hamiltonian model $\Kt$-spaces. First recall that for every $a \in \A$, the subgroup $K_{a\tau}$ has $\exp(\fa)$ as a maximal torus, whose character goup is $\lat_{\tau} \subset \fa$. The weight lattice of the complexification $K_{a\tau}^{\C}$ of $K_{a\tau}$ can naturally be identified with the weight lattice $\lat_{\tau}$ of $K_{a\tau}$. Furthermore $\R_{\ge 0}(\A-a) \subset \fa$ is a Weyl chamber for $K_{a\tau}$ and thus determines a Borel subgroup of $K_{a\tau}^{\C}$ with respect to which $\R_{\ge 0}(\A-a) \cap \lat_{\tau}$ is the set of dominant weights. In other words, when $X$ is a smooth affine spherical  $K_{a\tau}^{\C}$-variety, we view its weight monoid $\wm(X)$ as a subset of $\fa$.  

\begin{defn} \label{def:admissible}
Let $\fa = \ft^{\tau}$ and $\A \subset \fa$ be as in \cref{thm:twistconj}. Let $\lat$ be  a subgroup  of $\fa$. We will say that $\lat$ is \textbf{$\Kt$-admissible}  if for every vertex $a$ of $\A$  there exists a smooth affine spherical $K_{a\tau}^{\C}$-variety whose weight monoid $\Gamma_a$ satisfies $\Rp\Gamma_a = \Rp(\A - a)$ and $\Z\Gamma_a = \lat$. 
\end{defn}

\begin{rem}
\begin{enumerate}[(a),wide]
\item A subgroup $\lat$ of $\fa$ is $\Kt$-admissible if and only if, in the parlance of \cite[Definition 2.5.1]{knop-cmfqham-arxiv}, $(\A,\lat)$ is a \emph{spherical pair}.
\item If $\lat$ is $\Kt$-admissible, then the weight monoids $\Gamma_a$ as in \cref{def:admissible} are uniquely determined by $\lat$; see  \cref{rem:wm_adap_unique} below. Indeed, $\Gamma_a = \Rp(\A - a) \cap \lat$. 
\end{enumerate}
\end{rem}

Here is the announced specialization of  \cite[Corollary 2.6]{knop-cmfqham-arxiv}. 
\begin{thm}[Knop]\label{sphpair}
Let $\tau$ be an automorphism of the compact and simply connected Lie group $K$ and let $\fa = \ft^{\tau}$. 
The map $M \mapsto \lat_M$ yields a bijection between the set of isomorphism classes of quasi-Hamiltonian $\K\tau$-model spaces and the set of $\K\tau$-admissible subgroups of $\fa$.  
\end{thm}

We need a bit more notation before we can state our main result. 
Let $S_{\tau}$ be the set of simple roots of $\Phi_{\tau}$ corresponding to the choice of alcove $\A$, that is, $\alpha \in \Phi_{\tau}$ belongs to $S_{\tau}$ if and only if the affine hyperplane $H_{\alpha}$ is a wall of $\A$ and $\alpha(a)\ge 0$ for all $a \in \A$.  For $K$ simple, we number the simple roots in $S_{\tau} = \{\alpha_0, \alpha_1, \ldots,\alpha_n\}$ as in the Dynkin diagram $\mathsf{X}_n^{(r)}$ in \cref{ars} corresponding to $\Phi_{\tau}$.  

Here is the main result of this paper. The proof will be given in \cref{seclat}. 
\begin{thm} \label{thm_main}
Let $K$ be a simple and simply connected compact Lie group and $\tau$  a smooth automorphism of $K$. Let $T, \fa, \lat_\tau, \Phi_{\tau}$ and $\A$ be as in \cref{thm:twistconj} and number the  simple roots $S_{\tau}$ of $\Phi_{\tau}$ as in \cref{ars}.   Finally, let $\operatorname{ord}(\overline{\tau})$ be the order of the image of $\tau$ in the group of outer automorphisms of $K$. 
 
Then the map $M\mapsto \lat_M$ gives a bijection between the set of isomorphism classes of quasi-Hamiltonian model $K\tau$-spaces and the subgroups $\lat$ of $\fa$ in the following table:

\newcounter{acases}
\begin{longtable}{c|p{4cm}|c|p{9cm}|}
&$K$ & $\operatorname{ord}(\overline{\tau})$ &$\Kt$-admissible subgroups $\lat$ of $\fa$. \\
\hline
\refstepcounter{acases} \label{doubles} \theacases &
\multicolumn{2}{|p{6cm}|}{any $(K,\tau)$ except
$K=\SU(2n+1)$ with $n\ge 1$ and $\operatorname{ord}(\overline{\tau})=2$}& $\{2\aa_1,2\aa_2,\ldots,2\aa_n\} \subset \lat \subset 2\lat_{\tau}$\\
\hline
\refstepcounter{acases} \label{suodd} \theacases &
$\SU(n+1), n\ge 2$ even & 1 & $\{\aa_1,\aa_2,\ldots,\aa_n\} \subset \lat \subset \lat_{\tau}$
\\
\hline
\refstepcounter{acases}\label{sueven} \theacases &
$\SU(n+1), n\ge 1$ odd & 1 & $\la \aa_2+\aa_3, \aa_3+\aa_4,\dots, \aa_{n-1}+\aa_n, e\w_{n-1}, r\w_{n-1}+\w_n \ra_{\Z}$\newline with $r,e \in \Z_{\ge 1}$, $e|\frac{n+1}{2}$, $0 \le r \le e-1$, where $\w_{n-1}, \w_n \in \fa$ are defined by $\la \w_k, \aav_j\ra = \delta_{kj}$ for all $k \in \{n-1,n\}, j \in \{1,2,\ldots,n\}.$\\
\hline
\refstepcounter{acases} \theacases  \label{sp2nwl}&
$\Sp(2n), n\ge 2$ & 1 & $\lat_\tau$\\
\hline
\refstepcounter{acases} \theacases \label{su5}&
$\SU(5)$& 2 & $\la \aa_1, \aa_2\ra_{\Z}$\\
\hline
\refstepcounter{acases} \theacases \label{su2np12} &
$\SU(2n+1), n\ge 1$ & 2 & $\lat_{\tau}$\\
\hline
\refstepcounter{acases} \theacases \label{su2np12_2}&
$\SU(2n+1), n \ge 1$ & 2 & $2\lat_{\tau}$\\
\hline
\refstepcounter{acases} \theacases \label{spineven}&
$\Spin(2n+2), n\ge 2$ & 2 & $\la \aa_1,\aa_2,\ldots,\aa_n\ra_{\Z}$\\
\hline
\refstepcounter{acases} \theacases \label{spinodd}&
$\Spin(2n+2)$, {$n\ge 3$ odd} & 2 & $\la \aa_1+\aa_2, \aa_2+\aa_3, \dots, \aa_{n-1}+\aa_n, 2\aa_n\ra_{\Z}$\\
\hline
\end{longtable}
\end{thm}

\begin{rem}  \label{rem:after_thm_main}
We keep the notations from \cref{thm_main} in this remark.
\begin{enumerate}[(a), wide] 
\item When $\operatorname{ord}(\overline{\tau}) = 1$, $\lat_{\tau}$ is simply the weight lattice $\Xi(T)$ of $K$ and $\{\aa_1, \aa_2 , \ldots, \aa_n\}$ is the set of simple roots of $K$. This claim about $\lat_{\tau}$ holds because $\aa_0^{\vee}$ is an integral linear combination of $\aa_1^\vee,\aa_2^\vee,\ldots,\aa_n^\vee$, which holds, for example, because $-\aa_0$ is the highest root in the root system $\overline{\Phi}(\fk,\ft)$ of $K$ and the coroots of the simple roots form a basis of the dual root system.   \label{rem_wl_at_0_untwisted}
\item More generally one can check for each irreducible affine root system in \cref{ars} that $\lat_{\tau}$ is the weight lattice of the root system $\overline{\Phi_{\tau}(\sv_0)}$ of $K_{\sv_0 \tau}$; see \cref{lem:sphercheck}\ref{sphercheck1bis}. %\label{rem_wl_at_0}
\item The lattices $\lat$ in cases (\ref{doubles}) and (\ref{suodd}) are in natural bijective correspondence with the subgroups of the (finite) quotient $\lat_{\tau}/\la\aa_1,\aa_2,\ldots,\aa_n\ra_\Z$. For each irreducible finite root system, this quotient group is given in \cite[Planches I-IX]{bourbaki-geadl47}.
\item The following cases in \cref{thm_main} had already been found in \cite[Theorem 11.4]{knop-qham-arxiv}, see also \cite[Proposition 2.7.2]{knop-cmfqham-arxiv}: (\ref{sp2nwl}), (\ref{su2np12}), (\ref{sueven}) with $e=1, r=0$ and (\ref{suodd}) with $k=1$.
\item As is well known, and can be read in \cite[Planche I]{bourbaki-geadl47}, the weights $\w_{n-1}$ and $\w_n$ used to describe the lattices $\lat$ in case (\ref{sueven}) of \cref{thm_main} can be expressed as (rational) linear combinations of the simple roots $\aa_1,\aa_2,\ldots,\aa_n$. We have not found an elegant basis of $\lat$ in terms of these simple roots. Note that for $e=1$, the lattice is equal to $\lat_\tau$. 
\end{enumerate}
\end{rem}

\section{$G$-adapted lattices}\label{secloc}
Let $K$ be simple and simply connected.
If $\lat$ is a subgroup of $\fa$ that is $\Kt$-admissible and $\mathsf{v}$ is a vertex of $\A$, then there exists a smooth affine spherical $\Kvtc$-variety whose weight monoid $\Gamma_a$ satisfies $\Z\Gamma_a = \lat$ and $\Rp\Gamma_a = \Rp(\A-\sv)$. The first ingredient in the proof in \cref{seclat} of our classification of quasi-Hamiltonian model spaces is \cref{fullrank}, which was obtained in \cite[Section 3]{ppvs} and provides, up to replacing $K_{\sv_0\tau}^\C$ by its simply connected covering group $G(\sv_0)$, all the lattices that satisfy the condition for being $\Kt$-admissible at the vertex $\sv_0$ of $\A$ corresponding to the node $\a_0$ in the Dynkin diagram of $\Phi_{\tau}$ (the vertex $\sv_0$ is defined in \cref{eq:vi}). It will then remain to check which of these lattices verify the condition at every vertex of $\A$. When we do this in \cref{seclat}, we will make use of \cref{lem:sphercheck} and \cref{pvsalg}.  The latter is a special case of the combinatorial characterization of the weight monoids of smooth affine spherical varieties due to G. Pezzini and the second author, see \cite{charwm}. This  section also contains the necessary preliminaries to state \cref{pvsalg}.

\emph{For the remainder of this section,} $G$ is a complex connected reductive group, $B$ a chosen Borel subgroup, $H$ a chosen maximal torus in $B$, $\wl := \Hom(B,\C^{\times}) \equiv \Hom(H,\C^{\times})$ the weight lattice of $G$, $\oS$ the set of simple roots of $G$ with respect to $B$ and $H$ and $\dw$ the subset of dominant weights in $\wl$ with respect to $B$. Whenever necessary, we number the simple roots $\aa_1,\aa_2,\ldots \in \oS$ and the fundamental weights $\w_1,\w_2,\ldots \in \wl$ as in \cite[Planches I--IX]{bourbaki-geadl47}. 

\begin{defn} \label{def:Gadap}
Let $\Xi$ be a subgroup of $\wl$. We say that  $\Xi$ is  \textbf{$G$-adapted}  
if there exists a smooth affine spherical $G$-variety whose weight monoid $\Gamma$ satisfies 
\begin{align}
 \Z\Gamma &= \Xi, \text{ and} \label{eq:Gadap_1} \\
\R_{\ge 0} \Gamma &= \R_{\ge 0}\dw \text{ in }\wl \otimes_{\Z} \R  \label{eq:Gadap_2}
\end{align}
\end{defn} 

\begin{rem} \label{rem:wm_adap_unique}
\begin{enumerate}[(a), wide]
\item Because a smooth affine spherical variety is normal, its weight monoid $\wm$ satisfies the equality $\Gamma = \Rp\Gamma \cap \Z\Gamma$ in $\wl\otimes_{\Z}\R$. This means in particular that if $\Xi$ is $G$-adapted then there is only one monoid $\Gamma$ for which \cref{eq:Gadap_1,eq:Gadap_2} hold, namely $\Gamma = \Xi \cap (\R_{\ge 0}\dw)$.  \label{rem:wm_adap_unique_item1}
\item Furthermore, thanks to a theorem of I. Losev's in \cite{losev-knopconj}, a  smooth affine spherical $G$-variety $X$ is uniquely determined by its weight monoid (up to $G$-equivariant isomorphism). 
\end{enumerate}
\end{rem}

Part \ref{sphercheck3} of the following lemma allows us to ``ignore'' the lattice $\lat_{\tau}$ when determining whether a subgroup of $\fa$ is $\Kt$-admissible, making it a purely local problem at every vertex of $\A$.
\begin{lem} \label{lem:sphercheck}
We make the same assumptions as in \cref{thm_main}. For each vertex $\sv$ of $\A$ we let $\wl_{\sv} \subset \fa$ be the weight lattice of the root system $\overline{\Phi_{\tau}(\sv)}$. Then the following hold:
\begin{enumerate}[(a)]
\item $\bigcap_{\sv} \wl_{\sv} = \lat_{\tau}$, where the intersection is over all vertices $\sv$ of $\A$; \label{sphercheck1}
\item With the numbering of the simple roots $S_\tau = \{\a_0,\a_1,\ldots,\a_n\}$ of $\Phi_\tau$
as in \cref{ars}, we define the vertex $\sv_k$ of $\A$ by
\begin{equation} \label{eq:vi}
\{\sv_k\} := \{a \in A : \a(a)=0 \text{ for all }\a \in S_{\tau} \setminus \{\a_k\}\}
\end{equation} for each $k \in \{0,1,\ldots,n\}$. Then $\wl_{\sv_0} \subset \wl_{\sv_k}$ for all $k \in \{0,1,\ldots,n\}$. In particular, $\lat_{\tau} = \wl_{\sv_0}$.
\label{sphercheck1bis}
\item $(\overline{\Phi_{\tau}(\sv)},\wl_{\sv})$ is the integral root system of the simply connected covering group of $\Kvtc$, which we will denote by $G(\sv)$; \label{sphercheck2}
\item A subgroup $\Xi$ of $\fa$ is $\Kt$-admissible if and only if $\Xi$ is $G(\sv)$-adapted for every vertex $\sv$ of $\A$. \label{sphercheck3} 
\end{enumerate}
\end{lem}
\begin{proof}
Assertion \ref{sphercheck1} is essentially a restatement of part \ref{twistconjwl} of \cref{thm:twistconj}. The first assertion in \ref{sphercheck1bis} will be repeated and proved in \cref{lem:dlconseq}\ref{dlconseq_wl}. That $\lat_\tau = \wl_{\sv_0}$ follows because every vertex of $\A$ is of the form $\sv_k$ for some $k \in \{0,1,\ldots,n\}$. 
Assertion \ref{sphercheck2} is a standard fact of Lie theory.  We come to assertion \ref{sphercheck3}. The ``only if'' statement holds because if $X$ is a smooth affine spherical $\Kvtc$-variety, then the action lifts to $G(\sv)$. The ``if'' statement is true because if $\Xi$ is $G(\sv)$-adapted at every vertex $\sv$ of $\A$, then $\Xi$ lies in $\wl_{\sv}$ for every $\sv$. By  \ref{sphercheck1} it then follows that $\Xi \subset \lat_{\tau}$, which implies that at each vertex $\sv$, the $G(\sv)$-action on the smooth affine spherical $G(\sv)$-variety associated to $\sv$ factors through $\Kvtc$. 
\end{proof}

Part \ref{lem:Gadap_item3} of the next lemma gives a different description of $G$-adapted lattices. We will say that a subgroup $\Xi$ of $\wl$ has \textbf{full rank} if $\rank(\Xi) = \rank(\wl)$. Furthermore, we recall that a submonoid $\Gamma$ of $\dw$ is called \textbf{$G$-saturated} if $\Z\Gamma \cap \dw = \Gamma$. 
\begin{lem} \label{lem:Gadap}
 \begin{enumerate}[(a)]
\item If $\Xi$ is a subgroup of $\wl$, then $\Xi \cap \dw = \Xi \cap \Rp\dw$ (as subsets of $\wl \otimes_{\Z}\R$).  \label{lem:Gadap_item1}
\item If $\Xi$ is a subgroup of $\wl$ of full rank, then $\Rp(\Xi \cap \dw) = \Rp\dw$.  \label{lem:Gadap_item1bis}
\item The map $\Xi \mapsto \Xi \cap \dw$ is a bijection  from the set
\[\{\Xi: \Xi \text{ is a subgroup of full rank of }\wl\}\]  to the set  \[\{\Gamma: \Gamma \text{ is a $G$-saturated submonoid of $\dw$ with $\Z\Gamma$ of full rank}\}\]
with inverse map $\Gamma \mapsto \Z\Gamma$.  \label{lem:Gadap_item2}
\item A subgroup $\Xi$ of $\wl$ is $G$-adapted if and only if $\Xi$ is of full rank and $\Xi \cap \dw$ is the weight monoid of a smooth affine spherical $G$-variety. \label{lem:Gadap_item3}
\end{enumerate}
\end{lem}
\begin{proof}
Assertion \ref{lem:Gadap_item1} follows from the well-known fact that $\Rp\dw  \cap P = \dw$. Assertion \ref{lem:Gadap_item1bis} holds because every extremal ray of the convex polyhedral cone $\Rp\dw$ contains an element of $\dw$ and, since $\Xi$ has finite index in $\wl$, also an element of $\Xi$. 
Part \ref{lem:Gadap_item2} is a consequence of the fact that when $\Xi$ is a subgroup of full rank of $P$, then 
\begin{equation} \label{eq:latGsat}
\Z(\Xi \cap \dw) = \Xi.
\end{equation}
\Cref{eq:latGsat} in turn can be shown with essentially the same proof as  \cite[Prop.\ 1.1(iii)]{oda-convex}. 
We turn to assertion \ref{lem:Gadap_item3} and begin with the ``only if'' statement. Suppose that $\Xi$ is $G$-adapted. Since $\dw$ spans $\wl \otimes_{\Z} \R$ as a vector space, it follows from \cref{eq:Gadap_2} that $\Xi$ has full rank. Furthermore, it follows from \cref{rem:wm_adap_unique}\ref{rem:wm_adap_unique_item1} that $\Xi \cap \dw$ is the weight monoid of a smooth affine spherical $G$-variety. The ``if'' statement holds by \cref{eq:latGsat} and assertion \ref{lem:Gadap_item1bis}. 
\end{proof}

\Cref{fullrank} below summarizes Propositions 3.7 and  3.16 of \cite{ppvs}. Note that these two propositions in \emph{loc.cit.} are stated in terms of $G$-saturated submonoids of $\dw$ of full rank and that parts \ref{lem:Gadap_item2} and \ref{lem:Gadap_item3} of \cref{lem:Gadap} show that this is just different terminology for the same objects. We'll  make use of the following notation:
\begin{equation} \label{eq:2osplus}
2\oS:=\{2\aa:\aa\in \oS\} \quad \text{ and }\quad 
\oS^+:=\{\aa+\bet: \aa,\bet\in \oS, \aa\ne \bet, \aa \not\perp \bet\}.
\end{equation}
\begin{prop} \label{fullrank}
Suppose $G$ is simply connected and simple. A sublattice $\Xi$ of the weight lattice $\wl$ of $G$ is $G$-adapted  if and only if one of the following holds
\begin{enumerate}[label={(AL\arabic*)}]
\item $2\oS \subset \Xi \subset 2\wl$, \label{al_doubles}
\item $G$ is of type $\mathsf{A}_n$ with $n\ge 1$, $n$ even and $\oS^+ \subset \Xi$; \label{al_An_even}
\item $G$ is of type $\mathsf{A}_n$ with $n\ge 1$, $n$ odd, $\oS^+ \subset \Xi$ and the odd coroots \(\aa_1^{\vee}|_{\Xi}, \aa_3^{\vee}|_{\Xi}, \ldots, \aa_n^{\vee}|_{\Xi}\) are part  of a basis of the dual lattice $\Xi^*:=\Hom_{\Z}(\Xi, \Z)$; \label{al_An_odd}
\item $G$ is of type $\mathsf{B}_n$ with $n \ge 2$  and
\(\Xi=\Z(S^+ \cup \{2\aa_n\}) ;\)\label{al_Bn_1}
\item $G$ is of type $\mathsf{B}_n$ with $n \ge 2$  and
\(\Xi = \la \aa_1,\aa_2,\ldots,\aa_n\ra_{\Z};\)\label{al_Bn_2}
\item $G$ is of type $\mathsf{C}_n$ with $n\ge 2$ and $\Xi = \wl$. \label{al_Cn}
\end{enumerate}
\end{prop}

\begin{rem} \label{rem_after_fullrank}
\begin{enumerate}[(a),wide]
\item \label{rem_after_fullrank_B2C2} If $G$ is of type $\mathsf{B}_2 \cong \mathsf{C}_2$, then there are five $G$-adapted subgroups of $\wl$:
\begin{equation} \label{eq:fullrank_B2C2}
\wl,\quad 2\Z\oS,\quad 2\wl,\quad \la \aa_1+\aa_2, 2\aa_2 \ra_\Z , \quad \la\aa_1,\aa_2\ra_\Z,
\end{equation} 
where $\alpha_1$ is the long simple root and $\alpha_2$ the short one. 
\item The lattices $\Xi$ as in \ref{al_doubles} are in natural correspondence with the subgroups of the (finite) quotient $2\wl/2\Z\oS \cong \wl/\Z\oS$. For each simple and simply connected $G$, the group $\wl/\Z\oS$ is given in \cite[Planches I-IX]{bourbaki-geadl47}.
\item  \label{rem_after_fullrank_Anmod} For concreteness and as we will make use of it in what follows, we recall from \cite[Lemma 3.10]{ppvs} the explicit description   of the lattices in \ref{al_An_even} and \ref{al_An_odd}. Let $G$ be of type $\mathsf{A}_n$ with $n\ge 1$ and $\Xi$ a subgroup of $\wl$. 
\begin{itemize}
\item Suppose $n$ is even. Then $\Xi$  satisfies \ref{al_An_even} if and only if $\Xi = \Z\oS^+ \oplus \Z(k\w_{n-1})$ for some $k \in \N\setminus\{0\}$.
\item Suppose $n$ is odd. Then $\Xi$ satsifies \ref{al_An_odd} if and only if 
\[\Xi = \la \aa_2+\aa_3, \aa_3+\aa_4, \ldots,\aa_{n-1}+\aa_n, e\w_{n-1},r\w_{n-1} + \w_n\ra_{\Z},\]
for some $e,r \in \N$ with $e|\frac{n+1}{2}$ and $0\le r\le e-1$.
\end{itemize}
\item For each lattice $\Xi$ in \cref{fullrank}, Tables 2 and 3 in \cite{ppvs} contain an explicit description of the smooth affine spherical $G$-variety $X$ such that $\Z\G(X) = \Xi$. These provide ``local models'' of quasi-Hamiltonian model spaces, in the following sense.  Suppose $\lat$ is $\Kt$-admissible, let $(M,m)$ be the $\Kt$-model space determined by $\lat$ and let $a$ be a vertex of $\A$. If $X$ is the smooth affine spherical $K_{a\tau}^\C$-variety whose weight monoid is $\Rp(\A-a)\cap \lat$, then Remark 2.5.4(d) of \cite{knop-cmfqham-arxiv} explains how $X$ describes a neighborhood of $m_+^{-1}(a)$ in $M$.
\item \label{rem_after_fullrank_at_v} In \cref{seclat} we will use an expression like ``$\Xi$ satisfies \ref{al_doubles} at [the vertex] $\sv$ [of $\A$]'' to say that $\Xi$ is a $G(\sv)$-adapted subgroup of $\wl_{\sv}$ satisfying \ref{al_doubles} for $G=G(\sv)$. 
\end{enumerate}
\end{rem}

\Cref{fullrank} was proved in \cite{ppvs} by using the combinatorial characterization of the weight monoids of smooth affine spherical varieties from \cite{charwm}. We now present, in \cref{pvsalg}, a special case of Theorem 1.12 of \emph{loc.cit.}, which we will use in \cref{seclat} when we verify whether a lattice is $G(\sv)$-adapted for a group $G(\sv)$ which is not simple. We first need to introduce two objects.

\begin{defn}\label{nsph}
Let $\Xi$ be a sublattice of $\wl$ of full rank. We define the set $\SN(\Xi)$ of \textbf{N-spherical roots of $\Xi$} as follows:
\[
\SN(\Xi) = (\oS^+ \cap \Xi) \cup \{2\aa \in 2\oS\cap \Xi: \la \aa^\vee, \gamma \ra \in 2\Z\text{ for all }\gamma \in \Xi\},
\]
where $\oS^+$ and $2\oS$ are the sets defined in \cref{eq:2osplus}
\end{defn}

\begin{prop}[{see \cite[Prop 1.7]{charwm}}]  \label{prop:SXi}
Let $\Xi$ be a sublattice of $\wl$ of full rank. Among all the subsets $F\subseteq \oS$ such that the relative interior of the cone spanned by $\{\aav:\aa\in F\}$ in $\operatorname{Hom}_\Z(\Xi, \Q) = \Hom_{\Z}(\wl,\Q)$ contains a point $x$ with $\la \sigma,x\ra \le 0$ for all $\sigma \in \SN(\Xi)$ there is a unique one, denoted $\oS_\Xi$, which contains all the others.
\end{prop}

Here is the announced specialization of \cite[Theorem 1.12]{charwm}. In order to save space, we freely use notions from \cite[\S 1]{charwm} and \cite[\S\S 2 and 3]{ppvs} in its proof. For the convenience of the reader, we point out that everything we need from \cite[\S 1]{charwm} is also contained in \cite[\S 2]{ppvs}. 
\begin{prop} \label{pvsalg} 
 Let $\Xi$ sublattice of $\wl$ of full rank. Then $\Xi$ is $G$-adapted if and only if
 \begin{enumerate}[(1)]
  \item \label{cond_basis} $\{\aa^\vee|_{\Xi}:\aa\in \oS_{\Xi}\}$ is a subset of a basis of the dual lattice $\Xi^*$, 
  \item \label{cond_perp} if $\aa,\bet$ in $\oS_{\Xi}$ and $\aa \neq \bet$, then $\aa \perp \bet$, and 
  \item \label{cond_no_double} if $\aa \in \oS_{\Xi}$, then $2\aa \not\in \Sigma^N(\Xi)$. 
 \end{enumerate}
\end{prop}

\begin{proof}
By \cref{lem:Gadap}\ref{lem:Gadap_item3} it suffices to show that $\Xi$ satisfies the conditions \ref{cond_basis}, \ref{cond_perp}, \ref{cond_no_double} of the proposition if and only if the $G$-saturated monoid $\Gamma:=\Xi \cap \dw$ satisfies the conditions (a), (b), (c) of \cite[Theorem 1.12]{charwm}. 

We first show that the set $\SN(\Xi)$ of \cref{nsph} is the same set as $\SN(\Gamma)$ in \cite[Theorem 1.12]{charwm} and in \cite[\S 3]{ppvs}. By \cite[Lemma 3.2(b)]{ppvs}, $\SN(\Gamma) \subset \oS^+ \cup 2\oS$.  Furthermore, by \cite[Lemma 3.2(a)]{ppvs}, the set $\oS^p(\Gamma):=\{\aa \in \oS: \la\aa^\vee, \lambda\ra =0 \text{ for all }\lambda \in \Gamma\}$ is empty. It is now straightforward to check, using \cite[Prop.\ 1.7]{charwm}, that $\SN(\Gamma) = \SN(\Xi)$. 

The equality $\SN(\Gamma) = \SN(\Xi)$ immediately gives us that the set $\oS_{\Xi}$ in \cref{prop:SXi} is the same as the set $\oS_{\Gamma}$ in \cite[Prop.\ 1.7]{charwm}. Because $\Z\Gamma = \Xi$ by \cref{eq:latGsat}, it follows that condition  \ref{cond_basis} of the current proposition is the same as condition (a) in \cite[Theorem 1.12]{charwm}. 

Suppose now that $\G$ fulfills conditions (a), (b) and (c) of \cite[Theorem 1.12]{charwm}. Then, $\Xi$ satisfies \ref{cond_perp} and \ref{cond_no_double} of the current proposition by \cite[lemma 3.4.]{ppvs}. 

Conversely, suppose that $\Xi$ fulfills  \ref{cond_basis}, \ref{cond_perp}, \ref{cond_no_double}  of the current proposition.  Because $\Z\Gamma$ is of full rank, there are no simple roots $\aa,\bet \in \oS$ such that $\aa\neq \bet$ and $\aa^{\vee}|_{\Z\Gamma} = \bet^\vee|_{\Z\Gamma}$, and so condition (b) of \cite[Theorem 1.12]{charwm}  is trivially met. Finally, it follows from \ref{cond_no_double} that $\oS_{\Gamma} \cap \SN(\Gamma) = \oS_{\Xi} \cap \SN(\Xi) = \emptyset$. Together with \ref{cond_perp} this implies that the triple
$(\oS_{\Gamma}, \oS^p(\Gamma), \oS_{\Gamma} \cap \SN(\Gamma)) = (\oS_{\Gamma}, \emptyset, \emptyset)$ is a (possibly empty) ``disjoint union'' of copies of the triple $(\mathsf{A}_1, \emptyset, \emptyset)$ in \cite[List 1.10]{charwm}. In particular, the triple satisfies condition (c) of \cite[Theorem 1.12]{charwm}, and we have shown that $\Gamma$ satisfies all three conditions in \emph{loc.cit.}
\end{proof}

We conclude this section with a generalization of \ref{al_doubles}, which has the same proof as \cite[Prop. 3.7]{ppvs}
\begin{prop} \label{prop:doubles}
If $\Xi$ is a subgroup of full rank of $\wl$ satisfying $2\oS \subset \Xi \subset 2\wl$, then $\Xi$ is $G$-adapted.  
\end{prop}
\begin{proof}
Because $2\oS \subset \Xi \subset 2\wl$, we have $\SN(\Xi) = 2\oS$. One then computes, that $\oS_{\Xi} = \emptyset$. Consequenlty, the conditions in \cref{pvsalg} are trivially satisfied. 
\end{proof}

\section{Proof of \cref{thm_main}} \label{seclat}

In this section we will prove \cref{thm_main}. \emph{For the remainder of this paper,} $K$ is a simply connected and compact Lie group, $\tau$ an automorphism of $K$, and we fix $T, \fa, \lat_\tau, \Phi_{\tau}$ and $\A$ as in \cref{thm:twistconj}. As before we will use $S_{\tau}$ for the set of simple roots of $\Phi_{\tau}$ determined by the choice of alcove $\A$ and we will number the elements of $S_{\tau} = \{\a_0,\a_1,\a_2,\ldots,\a_n\}$ as in \cref{ars}. We will also use the notations $\wl_\sv$ and $G(\sv)$ from \cref{lem:sphercheck} and we set
\[\oS_{\tau}:=\{\aa: \a \in S_{\tau}\}.\]
If $\sv$ is a vertex of $\A$, we set
\[%\label{def:oSv}
\oS(\sv):=\{\aa \colon \a \in S_{\tau}, \a(\sv)=0\}.
\] 
Then $\oS(\sv)$ is the set of simple roots of $\Kvt$, $\Kvtc$ and $G(\sv)$ corresponding to the choice of $\Rp(\A - \sv)$ as the positive Weyl chamber.  
Let $\ell \in \{0,1,2,\ldots,n\}$. We recall the definition of the vertex  $\sv_\ell$ of $\A$ from \cref{eq:vi}. Then the Dynkin diagram of $\oS(\sv_\ell)$ is obtained by removing from the Dynkin diagram $\mathsf{X}_n^{(r)}$ of $S_{\tau}$ the simple root $\alpha_\ell$ and all the edges adjacent to it.  
The following notation will also be useful:
\[\oS(\sv_{\ell})^+:=\{\aa+\bet: \aa,\bet\in \oS(\sv_\ell), \aa\ne \bet, \aa \not\perp \bet\}.\]
Furthermore, we will use $\w_1,\w_2,\ldots,\w_n$ for the fundamental weights of the root system $\overline{\Phi_{\tau}(\sv_0)}$, that is, $(\w_i)_{i=1}^n$ are those elements of $\fa$ such that $\la \w_i, \aa^\vee_j\ra = \delta_{ij}$ for all $i,j \in \{1,2,\ldots,n\}$.
We will make frequent use of the expression ``satisfies \ref{al_doubles} at $\sv$'' introduced in \cref{rem_after_fullrank}\ref{rem_after_fullrank_at_v}.

We begin by explaining the Dynkin labels $k(\alpha)$ attached to the simple roots $\alpha$ in each diagram in \cref{ars}. They are the unique coprime positive integers such that 
\(\delta:= \sum_{\alpha \in S_{\tau}}k(\alpha)\alpha\)
is a constant function. Taking gradients we obtain the equation
\begin{equation} \label{eq:delta_grad}
\sum_{\alpha \in S_{\tau}} k(\alpha)\aa = 0
\end{equation}
which will be important in what follows. One immediate consequence, using the definition \eqref{eq:coroot} of $\aa^{\vee}$, is
\begin{equation} \label{eq:corootequation}
\sum_{\alpha \in S_{\tau}} k(\alpha)\|\aa\|^2\aa^\vee = 0.
\end{equation}
Since it will play a role, we recall that the number of edges between two simple roots in a Dynkin diagram  gives information about their relative lengths:
\begin{align*}
&\dynkin[labels={\alpha,\beta}, edge length=.7cm] A2 \quad \text{means} \quad  \frac{\|\aa\|^2}{\|\bet\|^2} = 1; 
&\dynkin[labels={\alpha,\beta}, edge length=.7cm] B2 \quad \text{means} \quad  \frac{\|\aa\|^2}{\|\bet\|^2} = 2; \\
&\dynkin[labels={\alpha,\beta}, edge length=.7cm] G2 \quad \text{means} \quad  \frac{\|\aa\|^2}{\|\bet\|^2} = 3; 
&\dynkin[labels={\beta, \alpha}, edge length=.7cm, backwards] A[2]2 \quad \text{means} \quad  \frac{\|\aa\|^2}{\|\bet\|^2} = 4. 
\end{align*}

The following summarizes some immediate consequences of \cref{eq:delta_grad,eq:corootequation}. 

\begin{lem} \label{lem:dlconseq}
\begin{enumerate}[(a)]
\item $\wl_{\sv_0} \subset \wl_{\sv_\ell}$ for all $\ell \in \{0,1,\ldots,n\}$.  \label{dlconseq_wl}
\item If $\ell \in \{0,1,,\ldots,n\}$ with $k(\a_{\ell})=1$ then $\Z\oS(\sv_{\ell}) = \Z\oS_\tau$. \label{dlconseq_rl_eq}
\item If $\ell \in \{0,1,\ldots,n\}$ with $k(\a_{\ell})>1$ then $\Z\oS(\sv_{\ell}) \subsetneq \Z\oS_\tau$. \label{dlconseq_rl_neq}
\end{enumerate}
\end{lem}
\begin{proof}
To show \ref{dlconseq_wl}  one uses \cref{eq:corootequation}  to check for each affine root system in \cref{ars} that, with the chosen numbering of the simple roots, $\aa_0^\vee$  is an integral linear combination of $\aa_1^\vee,\aa_2^\vee,\ldots,\aa_n^\vee$  (for the untwisted diagrams $\mathsf{X}_n^{(1)}$ one can also argue as in part \ref{rem_wl_at_0_untwisted} of \cref{rem:after_thm_main}.)
To show \ref{dlconseq_rl_eq} we first observe that $\oS(\sv_{\ell}) =\oS_{\tau} \setminus\{\aa_{\ell}\}$. The claim now follows because \cref{eq:delta_grad} implies that $\aa_{\ell} \in \Z\oS(\sv_\ell)$ when $k(\a_{\ell}) = 1$. Part \ref{dlconseq_rl_neq} follows from the fact that the Dynkin labels are coprime, which together with the linear independence of the simple roots in $\oS(\sv_{\ell})$ implies that $\aa_{\ell} \notin \Z\oS(\sv_{\ell})$ when $k(\a_\ell)>1$. 
\end{proof}

We now start the actual \emph{proof of \cref{thm_main}.} For each irreducible affine root system in \cref{ars} we will check which of the $G(\sv_0)$-adapted sublattices $\Xi$ of $\wl_{\sv_0}$ are $\Kt$-admissible. It is \cref{fullrank} which provides those $G(\sv_0)$-adapted sublattices. 

This next proposition will determine all the $\Kt$-admissible subgroups of $\fa$ for many of the root systems $\Phi_\tau$ and justify entry (\ref{doubles}) in \cref{thm_main}. 

\begin{lem} \label{lem:aldoubles_adap}
If  $\Phi_\tau$ is not of type $\mathsf{A}_{2n}^{(2)}$, with $n\ge 1$ and $\Xi$ is a subgroup of $\wl_{\sv_0}$ with $2\oS(\sv_0) \subset \Xi \subset 2\wl_{\sv_0}$ then $\Xi$ is $\Kt$-admisible.
\end{lem}
\begin{proof}
One checks in \cref{ars} that $k(\a_0)=1$, and so it follows  from \cref{lem:dlconseq} that
\[2\oS(\sv_\ell) \subset \Z(2\oS(\sv_0)) \subset \Xi \subset 2\wl_{\sv_0} \subset 2\wl_{\sv_\ell}\]
for all $\ell \in \{1,2,\ldots,n\}$. \Cref{prop:doubles} tells us that $\Xi$ is $G(\sv_\ell)$-adapted for all $\ell$. By \cref{lem:sphercheck}\ref{sphercheck3} we obtain that $\Xi$ is $\Kt$-admissible. 
\end{proof}

\subsection*{Case: $\Phi_\tau$ has type \( \mathsf{D}_n^{(1)} \text{ with $n\ge 4$}, \mathsf{E}_6^{(1)},\mathsf{E}_7^{(1)}, \mathsf{E}_8^{(1)},\ \mathsf{F}_4^{(1)}, \mathsf{G}_2^{(1)},\ \mathsf{E}_6^{(2)}, \text{ or } \mathsf{D}_4^{(3)}\)}

Because the only $G(\sv_0)$-adapted subgroups of $\wl_{\sv_0}$ for these affine root systems are those satisfying \ref{al_doubles} at $\sv_0$, \cref{lem:aldoubles_adap} yields the  following.  
\begin{cor} \label{cor:uic}
Suppose $\Phi_\tau$ has one of the following Dynkin types:
\[\mathsf{D}_n^{(1)} \text{ with $n\ge 4$},\ \mathsf{E}_6^{(1)}, \mathsf{E}_7^{(1)},\ \mathsf{E}_8^{(1)},\ \mathsf{F}_4^{(1)},\ \ \mathsf{G}_2^{(1)},\ \mathsf{E}_6^{(2)},\ \mathsf{D}_4^{(3)}.\] 
Then the $\Kt$-admissible subgroups $\Xi$ of $\fa$  are the sublattices $\Xi$ of $\wl_{\sv_0}$ with
$2\oS(\sv_0) \subset \Xi \subset \lat_{\tau}.$ 
\end{cor}

This shows that \cref{thm_main} contains all the $\Kt$-admissible lattices when $\Phi_{\tau}$ is of one of the types listed in \cref{cor:uic}.

\subsection*{Case: $\Phi_\tau$ is of type $\mathsf{A}_1^{(1)}$}
Here $G(\sv_0) \cong \SL(2)$ and the only $G(\sv_0)$-adapted lattice not satisfying \ref{al_doubles}
is $\wl_{\sv_0} = \lat_{\tau}$. Because $k(\aa_0)=k(\aa_1)=1$ and $\|\aa_0\|^2 = \|\aa_1\|^2$ it follows from \cref{eq:corootequation} that $\wl_{\sv_0} = \wl_{\sv_1}$, which is $G(\sv_1)$-admissible, because  $G(\sv_1) \cong \SL(2)$

Together with \cref{lem:aldoubles_adap} we have shown
\begin{lem}
If $\Phi_{\tau}$ is of type $\mathsf{A}_1^{(1)}$, then the $\Kt$-admissible subgroups of $\fa$ are
$\lat_{\tau}, 2\lat_{\tau}$ and  $4\lat_{\tau}$.
\end{lem}
We have justified entry (\ref{sueven}) for $n=1$ in \cref{thm_main} and shown that \cref{thm_main} contains all $\Kt$-admissible lattices when $\Phi_\tau$ is of type $\mathsf{A}_1^{(1)}$
\begin{rem}
The $\Kt$-admissible lattices for $K\cong \SU(2)$ (and the corresponding manifolds) are already contained in \cite[\S 11, Example 2]{knop-qham-arxiv}, see also \cite[\S 2.7]{knop-cmfqham-arxiv}. 
\end{rem}

\subsection*{Case: $\Phi_\tau$ is of type $\mathsf{A}_n^{(1)}$ with $n\ge 2$ even} 

\begin{lem}  \label{lem:An1even}
Suppose $\Phi_\tau$ is of type $\mathsf{A}_n^{(1)}$ with $n\ge 2$ even and let $\Xi$ be a $G(\sv_0)$-adapted sublattice of $\wl_{\sv_0}$ that does not satisfy \ref{al_doubles} at $\sv_0$. Then the following are equivalent
\begin{enumerate}[(a)]
\item \label{An1even_a} $\Xi$ is $\Kt$-admissible; 
\item \label{An1even_b} $\aa_0+\aa_1, \aa_n+\aa_0 \in \Xi$;
\item \label{An1even_c} $\oS(\sv_0) \subset \Xi$;
\item \label{An1even_d} $\oS_\tau \subset \Xi$;
\end{enumerate}
\end{lem}
\begin{proof}
It follows from the assumptions on $\Xi$ and from \cref{fullrank}, that $\Xi$ satisfies \ref{al_An_even} at $\sv_0$, i.e.\ $\oS(\sv_0)^+ \subset \Xi$. 
We first show that \ref{An1even_a} and \ref{An1even_b} are equivalent.
Let $\ell \in \{0,1,\ldots,n\}$. The Dynkin type of $\oS(\sv_\ell)$ is $\mathsf{A}_n$  and $\oS(\sv_0)^+ \not\subset 2\wl_{\sv_\ell}$, because $(\aa_1^\vee, \aa_1+\aa_2) = (\aa_2^\vee,\aa_1+\aa_2) = 1$, which is odd. This implies that $\Xi \not \subset 2\wl_{\sv_\ell}$ and  it follows, again from \cref{fullrank}, that $\Xi$ is $G(\sv_{\ell})$-adapted if and only if $\oS(\sv_{\ell})^+ \subset \Xi$.  The equivalence of \ref{An1even_a} and \ref{An1even_b} now follows, with \cref{lem:sphercheck}\ref{sphercheck3}, from the  fact that  
$\bigcup_{\ell=0}^n \oS(\sv_{\ell})^+ \setminus \oS(\sv_{0})^+ = \{\aa_0+\aa_1, \aa_n+\aa_0\}$.

We now show that \ref{An1even_c} implies \ref{An1even_d}.
For this root system, \cref{eq:delta_grad} becomes
\begin{equation} \label{eq:delta_grad_A1n}
\aa_0 + \aa_1 + \cdots+\aa_n = 0,
\end{equation}
which implies, since $\Xi$ is a subgroup of $\wl_{\sv_0}$, that if $\Xi$ contains $\oS(\sv_0) = \{\aa_1,\aa_2,\ldots,\aa_n\}$, then it also contains $\oS_{\tau} = \oS(\sv_0) \cup \{\aa_0\}$. 
That \ref{An1even_b} follows from \ref{An1even_d} is clear. Finally, we prove that \ref{An1even_b} implies \ref{An1even_c}.  Observe that, since $n$ is even, \cref{eq:delta_grad_A1n} can be rewritten as
\[-\aa_n=\sum_{k=1}^{n/2}(\aa_{2(k-1)}+\aa_{2k-1}) = (\aa_0+\aa_1)+(\aa_2+\aa_3)+\cdots+(\aa_{n-2}+\aa_{n-1}).\]
If we now assume that \ref{An1even_b} holds, and in particular that $\aa_0+\aa_1 \in \Xi$, then this equation implies that $\aa_n \in \Xi$, since $\oS(\sv_0)^+ \subset \Xi$. Again using that $\oS(\sv_0)^+ \subset \Xi$ one then (recursively) deduces \ref{An1even_c}. 
\end{proof}

\begin{rem}
It follows from straightforward computations like in the proof of \cref{lem:An1even} that, under the assumptions of the lemma, assertion \ref{An1even_b} of the lemma holds if and only if $\aa_0+\aa_1 \in \Xi$ if and only if $\aa_n+\aa_0 \in \Xi$. 
\end{rem}

Together with \cref{lem:aldoubles_adap} we have proven
\begin{lem}
Suppose $\Phi_\tau$ is of type $\mathsf{A}_n^{(1)}$ with $n\ge 2$ even. The subgroups of $\fa$ that are $\Kt$-admissible are the lattices $\Xi$ with $2\oS(\sv_0) \subset \Xi \subset 2\lat_{\tau}$ and the lattices $\Xi$ with $\oS(\sv_0) \subset \Xi \subset \lat_{\tau}$. 
\end{lem}

We have justified entry (\ref{suodd}) in \cref{thm_main} and shown that \cref{thm_main} contains all $\Kt$-admissible lattices when $\Phi_\tau$ is of type $\mathsf{A}_n^{(1)}$ with $n\ge 2$ even.

\subsection*{Case: $\Phi_\tau$ is of type $\mathsf{A}_n^{(1)}$ with $n\ge 3$ odd}

\begin{lem}  \label{lem:An1odd}
Suppose $\Phi_\tau$ is of type $\mathsf{A}_n^{(1)}$ with $n\ge 3$ odd and let $\Xi$ be a $G(\sv_0)$-adapted sublattice of $\wl_{\sv_0}$ that does not satisfy \ref{al_doubles} at $\sv_0$. Then $\Xi$ is $\Kt$-admissible if and only if 
$\aa_0+\aa_1, \aa_n+\aa_0 \in \Xi$ and the even coroots $\aa_0^\vee, \aa_2^\vee,\ldots,\aa_{n-1}^\vee$ are part of a $\Z$-basis of $\Xi^*$. 
\end{lem}
\begin{proof}
As $\Xi$ is $G(\sv_0)$-adapted and does not satisfy \ref{al_doubles}, \cref{fullrank} implies that $\Xi$ satisfies \ref{al_An_odd} at $\sv_0$. In particular, it contains $\aa_1+\aa_2$. Let $\ell \in \{0,1,\ldots,n\}$. The Dynkin type of $\oS(\sv_\ell)$ is $\mathsf{A}_n$ and $\aa_1 \oplus \aa_2 \notin 2\wl_{\sv_\ell}$ because $\la \aa_1^\vee, \aa_1+\aa_2\ra = \la \aa_2^\vee, \aa_1+\aa_2\ra =1$. \Cref{fullrank} now yields that $\Xi$ is $G(\sv_\ell)$-adapted if and only if $\Xi$ satisfies \ref{al_An_odd} at $\sv_{\ell}$. 

With \cref{lem:Gadap}\ref{lem:Gadap_item3} and the fact that $\bigcup_{\ell=0}^n \oS(\sv_{\ell})^+ \setminus \oS(\sv_{0})^+ = \{\aa_0+\aa_1, \aa_n+\aa_0\}$, the lemma now follows because the coroots of $G=G(\sv_\ell)$ listed in \ref{al_An_odd} are $\aa_0^\vee,\aa_2^\vee, \ldots,\aa_{n-1}^\vee$ when $\ell$ is odd, and  $\aa_1^\vee,\aa_3^\vee, \ldots,\aa_{n}^\vee$ when $\ell$ is even. 
\end{proof}

The next lemma says that in this case all $G(\sv_0)$-adapted lattices which do not satisfy \ref{al_doubles} at $\sv_0$ are $\Kt$-admissible.
\begin{lem}
Suppose $\Phi_\tau$ is of type $\mathsf{A}_n^{(1)}$ with $n\ge 3$ odd and let $\Xi$ be a $G(\sv_0)$-adapted sublattice of $\wl_{\sv_0}$ that does not satisfy \ref{al_doubles} at $\sv_0$. Then
\begin{multline}\label{eq:al3latt}
\Xi = \la \aa_2+\aa_3, \aa_3+\aa_4, \ldots,\aa_{n-1}+\aa_n, e\w_{n-1},r\w_{n-1} + \w_n\ra_{\Z},\\ \text{ 
for some $e,r \in \N$ with $e|\frac{n+1}{2}$ and $0\le r\le e-1$}
\end{multline} 
and $\Xi$ is $\Kt$-admissible. 
\end{lem}
\begin{proof}
\Cref{eq:al3latt} follows from \cref{rem_after_fullrank}\ref{rem_after_fullrank_Anmod}.  By \cref{eq:delta_grad} we have
\begin{align*}
-(\aa_0+\aa_1)&=(\aa_2+\aa_3)+(\aa_4+\aa_5)+\cdots+(\aa_{n-1}+\aa_n)\text{ and }\\-(\aa_n+\aa_0)&=(\aa_1+\aa_2)+(\aa_3+\aa_4)+\dots+(\aa_{n-2}+\aa_{n-1}).
\end{align*}
Since $\oS(\sv_0)^+ \subset \Xi$, because $\Xi$ satsifies \ref{al_An_odd} at $\sv_0$, these two equations imply that $\aa_0+\aa_1, \aa_n+\aa_0 \in \Xi$. By \cref{lem:An1odd}, what remains is to show that the even coroots $\aa_0^\vee, \aa_2^\vee,\ldots,\aa_{n-1}^\vee$ are part of a $\Z$-basis of $\Xi^*$. To do so, we will apply the elementary divisors theorem, see e.g. \cite[Theorem 5.2, p. 234]{lang_algebra_second_ed}. We first recall that $\aa_0^{\vee} = -(\aa_1^{\vee}+\aa_2^{\vee}+\cdots +\aa_n^\vee)$ by \cref{eq:corootequation}.

Next we give the basis elements of $\Xi$ in \cref{eq:al3latt} a name, that is, for  $i \in \{1,2,\ldots,n\}$ we set
\[\sigma_i = \begin{cases}
\aa_{i+1} + \aa_{i+2} &\text{if $1\le i \le n-2$};\\
e\w_{n-1} &\text{if $i=n-1$};\\
r\w_{n-1}+\w_n &\text{if $i=n$}
\end{cases}\]
and we consider the matrix $A$ with  $n$ rows and   $d:=\frac{n+1}{2}$ columns and whose $(i,j)$-th entry is
\[
A_{ij} = \la\sigma_i, \aa^\vee_{2(j-1)}\ra.
\]
Put differently, the columns of $A$ give the coordinates of the coroots $\aa_0^{\vee}, \aa_2^\vee, \ldots,\aa_{n-1}^\vee$ in the basis of $\Xi^*$ that is dual to the basis $(\sigma_i)_{i=1}^n$ of $\Xi$. For example, for $n=7$ we have
\[A = \begin{pmatrix}
0 & 1 & -1  & 0 \\
0 & -1 & 1 & 0 \\
0 & 0 & 1 & -1 \\
0 & 0 & -1 & 1 \\
-1 & 0 & 0 & 1 \\
-e & 0 & 0 & e\\
-r-1 & 0 & 0 &r
\end{pmatrix}.
\]
We need to show that the greatest common divisor of all $d\times d$-minors of $A$ is $1$. 
To do so, we consider the $d\times d$-submatrix $M$ of $A$ consisting of rows $1,3,5, \ldots,n-2$ and $n$ of $A$.  For example for $n=7$, we have 
\[M = \begin{pmatrix}
0 & 1 & -1 & 0 \\
0 & 0 & 1 & -1 \\
-1 & 0 & 0 & 1 \\
-r-1 & 0 & 0 & r
\end{pmatrix}.
\]
Elementary computations show that $\det(M) = \pm 1$ and then the elementary divisors theorem implies that the even coroots are part of a $\Z$-basis of $\Xi^*$.
\end{proof}

Together with \cref{lem:aldoubles_adap} we have proven
\begin{lem}
Suppose $\Phi_\tau$ is of type $\mathsf{A}_n^{(1)}$ with $n\ge 3$ odd. The subgroups of $\fa$ that are $\Kt$-admissible are: the lattices $\Xi$ with $2\oS(\sv_0) \subset \Xi \subset 2\lat_{\tau}$ and the lattices in \cref{eq:al3latt}. 
\end{lem}
We have justified entry (\ref{sueven}) for $n\ge 3$ in \cref{thm_main} and shown that \cref{thm_main} contains all $\Kt$-admissible lattices when $\Phi_\tau$ is of type $\mathsf{A}_n^{(1)}$ with $n\ge 3$ odd.

\subsection*{Case: $\Phi_\tau$ is of type $\mathsf{B}_n^{(1)}$ with $n\ge 3$}

\begin{lem} \label{lem:Bn1}
Suppose $\Phi_\tau$ is of type $\mathsf{B}_n^{(1)}$ with $n\ge 3$ and let $\Xi$ be a $G(\sv_0)$-adapted sublattice of $\wl_{\sv_0}$ that does not satisfy \ref{al_doubles} at $\sv_0$. Then $\Xi$ is not $G(\sv_n)$-adapted, and therefore not $\Kt$-admissible. 
\end{lem}
\begin{proof}
For this affine root system $\Phi_\tau$, $G(\sv_0)$ is of type $\mathsf{B}_n$, with $n\ge 3$. As $\Xi$ does not satisfy \ref{al_doubles} at $\sv_0$, it follows that it satisfies \ref{al_Bn_1} or \ref{al_Bn_2} at $\sv_0$ which implies that 
\begin{equation}
\Xi \not\subset 2\wl_{\sv_n}. \label{eq:Bn1_lat}
\end{equation}  Indeed, this holds in both cases because $\aa_1+\aa_2 \in \Xi$, $\aa_1 \in \oS(\sv_n)$ and $\la\aa_1^\vee, \aa_1+\aa_2\ra = 1$ is odd. 

If $n\ge 4$, then $G(\sv_n)$ is of type $\mathsf{D_n}$, which means, by \cref{fullrank}, that the only $G(\sv_n)$-adapted lattices are those satisfying \ref{al_doubles}. By \eqref{eq:Bn1_lat} it follows that $\Xi$ is not $G(\sv_n)$-adapted.

If $n=3$, then $G(\sv_n) = G(\sv_3)$ is of type $\mathsf{A}_3$ and its Dynkin diagram is
\[ \dynkin[labels={0,2,1}, label macro/.code={\aa_{\drlap#1}}, edge length=.7cm] A{3} \]
By \eqref{eq:Bn1_lat}, $\Xi$ does not satisfy \ref{al_doubles} at $\sv_3$.  We show that it also doesn't satisfy \ref{al_An_odd} at $\sv_3$, which then implies by \cref{fullrank} that $\Xi$ is not $G(\sv_3)$-adapted, as there are no other adapted lattices for a group of type $\mathsf{A}_3$. 

If $\Xi$ satisfies \ref{al_Bn_1} at $\sv_0$, then $\Xi = \la \aa_1+\aa_2, \aa_2+\aa_3, 2\aa_3\ra_{\Z}$. Since the greatest common divisor of the $2\times 2$-minors of the matrix
\[
\begin{pmatrix}
\la \aa_0^\vee, \aa_1+\aa_2\ra & \la \aa_1^\vee, \aa_1+\aa_2\ra \\
\la \aa_0^\vee, \aa_2+\aa_3\ra & \la \aa_1^\vee, \aa_2+\aa_3\ra \\
\la \aa_0^\vee, 2\aa_3\ra & \la \aa_1^\vee, 2\aa_3\ra \\
\end{pmatrix}
= 
\begin{pmatrix}
-1 & 1 \\
-1 & -1 \\
0 & 0
\end{pmatrix}
\]
is $2$, the elementary divisors theorem tells us, that the coroots $\aa_0^\vee$ and $\aa_1^\vee$ are not part of a basis of the dual lattice $\Xi^*$. Consequently, $\Xi$ does not satisfy \ref{al_An_odd} at $\sv_3$ in this case.

If $\Xi$ satisfies \ref{al_Bn_2} at $\sv_0$, then $\Xi = \la \aa_1,\aa_2,\aa_3\ra_{\Z}$. 
Using the Dynkin diagram of $\mathsf{B}_3^{(1)}$, one computes the matrix
\[
\begin{pmatrix}
\la \aa_0^\vee, \aa_1\ra & \la \aa_1^\vee, \aa_1\ra \\
\la \aa_0^\vee, \aa_2\ra & \la \aa_1^\vee, \aa_2\ra \\
\la \aa_0^\vee, \aa_3\ra & \la \aa_1^\vee, \aa_3\ra \\
\end{pmatrix}
= 
\begin{pmatrix}
0 & 2 \\
-1 & -1 \\
0 & 0
\end{pmatrix}.
\]
As the greatest common divisor of the $2\times 2$-minors of this matrix is $2$, it follows again that $\{\aav_0, \aav_1\}$ is not part of a basis of $\Xi^*$, so that once again $\Xi$ does not satisfy \ref{al_An_odd} at $\sv_3$. 
\end{proof}

\Cref{lem:Bn1} and \cref{lem:aldoubles_adap} establish the following
\begin{lem}
Suppose $\Phi_\tau$ is of type $\mathsf{B}_n^{(1)}$ with $n\ge 3$. The subgroups of $\fa$ that are $\Kt$-admissible are the lattices $\Xi$ with $2\oS(\sv_0) \subset \Xi \subset 2\lat_{\tau}$. 
\end{lem}

This shows that \cref{thm_main} contains all $\Kt$-admissible subgroups of $\fa$ when $\Phi_{\tau}$ is of type $\mathsf{B}_n^{(1)}$ with $n\ge 3$. 

\subsection*{Case: $\Phi_\tau$ is of type $\mathsf{C}_n^{(1)}$ with $n\ge 2$}

Here  $G(\sv_0)$ is of type $\mathsf{C}_n$. By \cref{fullrank} (and \cref{rem_after_fullrank}\ref{rem_after_fullrank_B2C2}), the $G(\sv_0)$-adapted lattices are 
\begin{align}
&2\Z\oS(\sv_0), 2\wl_{\sv_0} \text{ and }\wl_{\sv_0}=\lat_\tau \text{ for all $n\ge 2$, and} \label{eq:cn1_gen} \\
&\text{in addition }\la \aa_1+\aa_2, 2\aa_1\ra_{\Z}\ \text{ and }\ \la \aa_1 , \aa_2\ra_{\Z} \text{ when $n=2$}. \label{eq:cn1_2}
\end{align}
We first deal with the lattices in \eqref{eq:cn1_gen}. The first two, $2\Z\oS(\sv_0), 2\wl_{\sv_0} $ are $\Kt$-admissible by \cref{lem:aldoubles_adap} and
it was shown in \cite[Proposition 2.7.2]{knop-cmfqham-arxiv} that $\lat_{\tau}$ is $\Kt$-admissible. This shows
\begin{lem} \label{lem:cn1_n3}
If  $\Phi_\tau$ is of type $\mathsf{C}_{n}^{(1)}$ with $n \ge 3$, then the $\Kt$-admissible subgroups of $\fa$ are $\lat_{\tau}$, $2\lat_{\tau}$ and $2\Z\oS(\sv_0)$.
\end{lem}
This justifies the entry (\ref{sp2nwl}) for $n\ge 3$ in \cref{thm_main} and shows that  \cref{thm_main} contains all $\Kt$-admissible lattices when $\Phi_\tau$ is of type $\mathsf{C}^{(1)}_{n}$ with $n\ge 3.$

\begin{lem} \label{lem:c21}
If  $\Phi_\tau$ is of type $\mathsf{C}_{2}^{(1)}$, then the $\Kt$-admissible subgroups of $\fa$ are 
\[\lat_{\tau},  2\lat_{\tau} \text{ and }2\Z\oS(\sv_0). \]
\end{lem}
\begin{proof}
Because the argument before \cref{lem:cn1_n3} also applies to the case $n=2$ we only need to consider the two lattices in \cref{eq:cn1_2}. We show that neither of them is $\Kt$-admissible using \cref{lem:sphercheck}\ref{sphercheck3}. 
We first show that 
\(\Xi= \la \aa_1 + \aa_2, 2\aa_1\ra_{\Z}   = \la \w_2,4\w_1 - 2\w_2 \ra_{\Z} = \la 4\w_1,\w_2\ra_\Z\) is not $G(\sv_1)$-adapted, using \cref{pvsalg}. Note that $G(\sv_1)$ is of type $\mathsf{A}_1 \times \mathsf{A}_1$. Writing $\wt_0$ and $\wt_2$ for the fundamental weights of $\wl_{\sv_1}$ we compute, using \cref{eq:corootequation}, that $\Xi= \la -4\wt_0, -\wt_0 + \wt_2\ra_\Z$.  An easy computation shows $\SN(\Xi) = \emptyset$ and it follows that $\oS_\Xi = \{\aa_0,\aa_2\}$. 
Because
\[\det \begin{pmatrix}
\la\aav_0, -4\wt_0 \ra & \la\aav_2, -4\wt_0 \ra  \\
\la\aav_0, -\wt_0+\wt_2 \ra & \la\aav_2, -\wt_0+\wt_2 \ra 
\end{pmatrix} = 
\det \begin{pmatrix}
-4 & 0 \\
-1 & 1
\end{pmatrix} = -4
\]
$\{\aav_0, \aav_2\}$ is not a basis of $\Xi^*$, and therefore condition \ref{cond_basis} of \cref{pvsalg} is not satisfied. 

Next, we show that $\Xi=\la 2\w_1,\w_2\ra_{\Z}$ is not $G(\sv_1)$-adapted. Again writing $\wt_0, \wt_2$ for the fundamental weights of $\wl_{\sv_1}$ we find that \begin{equation}
\Xi = \la -2\wt_0,-\ww_0 + \wt_2\ra_{\Z}. \label{eq:ansl2sl2lat}
\end{equation}
 Here too $\SN(\Xi)=\emptyset$ and therefore $\oS_\Xi=\{\aa_0,\aa_2\}$. Because
\[\det \begin{pmatrix}
\la\aav_0, -2\wt_0 \ra & \la\aav_2, -2\wt_0 \ra  \\
\la\aav_0, -\wt_0+\wt_2 \ra & \la\aav_2, -\wt_0+\wt_2 \ra 
\end{pmatrix} = 
\det \begin{pmatrix}
-2 & 0 \\
-1 & 1
\end{pmatrix} = -2
\]
$\{\aav_0, \aav_2\}$ is not a basis of $\Xi^*$, and once again condition \ref{cond_basis} of \cref{pvsalg} is not satisfied. 
\end{proof}
We have justified entry (\ref{sp2nwl}) for $n=2$ in \cref{thm_main} and shown that \cref{thm_main} contains all $\Kt$-admissible subgroups of $\fa$ when $\Phi_\tau$ is of type $\mathsf{C}^{(1)}_{2}$.

\subsection*{Case: $\Phi_\tau$ is of type $\mathsf{A}^{(2)}_{2}$}
Here $G(\sv_0) \cong \SL(2)$ and it is  well known (or can be read from \cref{fullrank}) that the 
$G(\sv_0)$-adapted lattices are 
\begin{equation}\Z\w_1, 2\Z\w_1\text{ and }4\Z\w_1. \label{eq:sl2lats} \end{equation} 
For this affine root system, \cref{eq:corootequation} reads
\[ 2\aa^\vee_0 + 4\aa_1^\vee = 0\]
so that $-2\aa_1^\vee =\aa_0^\vee$. This implies that $(\aa_0^\vee,\w_1) = -2$, so that 
\[\w_1 = -2\wt_0\]
where $\wt_0$ is the fundamental weight of $G(\sv_1)$. Since $G(\sv_1) \cong \SL(2)$ the only $G(\sv_1)$-adapted lattices are  \begin{equation}\Z\wt_0, 2\Z\wt_0\text{ and }4\Z\wt_0. \end{equation} 
and so of the three lattices in \eqref{eq:sl2lats} only
\[\Z\w_1 = 2\Z\wt_0 \text{ and } 2\Z\w_1= 4\Z\wt_0\]
are $\Kt$-admissible. 

We have thus proved:
\begin{lem}
If  $\Phi_\tau$ is of type $\mathsf{A}_{2}^{(2)}$, then the $\Kt$-admissible subgroups of $\fa$ are $\lat_{\tau}$ and $2\lat_{\tau}$.
\end{lem}
This justifies the entries (\ref{su2np12}) and (\ref{su2np12_2}) for $n=1$ in \cref{thm_main} and shows that \cref{thm_main} contains all $\Kt$-admissible lattices when $\Phi_\tau$ is of type $\mathsf{A}^{(2)}_{2}$.

\begin{rem}
The $\Kt$-admissible lattices when $\Phi_{\tau}$ is of type $\mathsf{A}^{(2)}_2$ are already contained in \cite[\S 11, Example 2]{knop-qham-arxiv}, see also \cite[\S 2.7]{knop-cmfqham-arxiv}. 
\end{rem}

\subsection*{Case: $\Phi_\tau$ is of type $\mathsf{A}_{2n}^{(2)}$ with $n \ge 2$}
Here  $G(\sv_0)$ is of type $\mathsf{C}_n$. By \cref{fullrank} (and \cref{rem_after_fullrank}\ref{rem_after_fullrank_B2C2}), the $G(\sv_0)$-adapted lattices are 
\begin{align}
&2\Z\oS(\sv_0), 2\wl_{\sv_0} \text{ and }\wl_{\sv_0}=\lat_\tau \text{ for all $n\ge 2$, and} \label{eq:a2n2_gen} \\
&\text{in addition }\la \aa_1+\aa_2, 2\aa_1\ra_{\Z}\ \text{ and }\ \la \aa_1 , \aa_2\ra_{\Z} \text{ when $n=2$}. \label{eq:a2n2_2}
\end{align}
We first deal with the lattices in \eqref{eq:a2n2_gen}. 
It was shown in \cite[Proposition 2.7.2]{knop-cmfqham-arxiv} that $\lat_{\tau}$ is $\Kt$-admissible. Since $k(\a_n)=1$ and $k(\a_0)=2>1$,  \cref{lem:dlconseq}\ref{dlconseq_rl_neq} tells us that $2\Z\oS(\sv_0)$ is a proper sublattice of $2\Z\oS(\sv_n)$. As $\oS(\sv_n)$ is of type $\mathsf{B}_n$ and $\Xi = 2\Z\oS(\sv_0)$ does not satisfy \ref{al_doubles}, \ref{al_Bn_1} or \ref{al_Bn_2}, nor \ref{al_Cn} when $n=2$, for $G=G(\sv_n)$,   it follows from \cref{fullrank}, that $2\Z\oS(\sv_0)$ is not $G(\sv_n)$-adapted and therefore also not $\Kt$-admissible. This leaves the $G(\sv_0)$-adapted lattice $2\wl_{\sv_0}$. It follows from \cref{lem:dlconseq}\ref{dlconseq_wl} that $$2\wl_{\sv_0} \subset 2\wl(\sv_k)\quad \text{for all $k \in \{1,2,\ldots,n\}$}.$$
Since $(\aa_i^\vee, \aa_\ell) \in \Z$  for all $i,\ell \in \{0,1,2,\ldots,n\}$ by \cref{def_ars}\ref{def_ars_integral},  we have $(\aa_i^\vee, 2\aa_\ell) \in 2\Z$ and consequently that
$$2\oS(\sv_k) \subset 2\wl_{\sv_0}\quad \text{for all $k \in \{1,2,\ldots,n\}$}.$$
By \cref{prop:doubles} and \cref{lem:sphercheck}\ref{sphercheck3} we obtain that $2\wl_{\sv_0} = 2\lat_\tau$ is $\Kt$-admissible. We have proved:
\begin{lem} \label{lem:A2n3}
If  $\Phi_\tau$ is of type $\mathsf{A}_{2n}^{(2)}$ with $n \ge 3$, then the $\Kt$-admissible subgroups of $\fa$ are $\lat_{\tau}$ and $2\lat_{\tau}$.
\end{lem}
This justifies the entries (\ref{su2np12}) and (\ref{su2np12_2}) for $n\ge 3$ in \cref{thm_main} and shows that  \cref{thm_main} contains all $\Kt$-admissible lattices when $\Phi_\tau$ is of type $\mathsf{A}^{(2)}_{2n}$ with $n\ge 3.$

\begin{lem}
If  $\Phi_\tau$ is of type $\mathsf{A}_{4}^{(2)}$, then the $\Kt$-admissible subgroups of $\fa$ are 
\[\lat_{\tau}, \ 2\lat_{\tau} \text{ and }\la \aa_1,\aa_2\ra_{\Z} \]
\end{lem}
\begin{proof}
Because the argument before \cref{lem:A2n3} also applies to the case $n=2$ we only need to consider the two lattices in \cref{eq:a2n2_2}. We fist show that 
\( \la \aa_1 + \aa_2, 2\aa_1\ra_{\Z}   = \la \w_2,4\w_1 - 2\w_2 \ra_{\Z}\) is not $G(\sv_2)$-adapted. Indeed, writing $\wt_0$ and $\wt_1$ for the fundamental weights of $\wl_{\sv_2}$ we compute, using \cref{eq:corootequation}, that $ \la \aa_1 + \aa_2, 2\aa_1\ra_{\Z} = \la 4\wt_1, 2\wt_0\ra_\Z$.  As this lattice does not satisfy \ref{al_doubles}, \ref{al_Bn_1}, \ref{al_Bn_2} or \ref{al_Cn} at $\sv_2$, it is not $G(\sv_2)$-adapted. 

Next, we show that $\Xi=\la \aa_1,\aa_2\ra_\Z=\la 2\w_1,\w_2\ra_{\Z}$ is $\Kt$-admissible, by showing that it is $G(\sv_1)$- and $G(\sv_2)$-adapted. Writing $\wt_0, \wt_1$ for the fundamental weights of $\wl_{\sv_2}$ we compute, using  \cref{eq:corootequation}, that $\Xi = 2\wl_{\sv_2}$, which is a $G(\sv_2)$-adapted lattice by \ref{al_doubles}. To show that $\Xi$ is $G(\sv_1)$-adapted, we use \cref{pvsalg}. If we now write $\wt_0, \wt_2$ for the fundamental weights of $\wl_{\sv_1}$ we find that $\Xi = \la -4\wt_0,-2\ww_0 + \wt_2\ra_{\Z}$. Straightforward computations show that $\SN(\Xi)=\{2\aa_0\}$ and $\oS_\Xi=\{\aa_2\}$ and then also that the three conditions in \cref{pvsalg} for $\Xi$ to be $G(\sv_1)$-adapted are satisfied. 
\end{proof}
This justifies the entries (\ref{su5}), (\ref{su2np12}) and (\ref{su2np12_2}) for $n=2$ in \cref{thm_main} and shows that \cref{thm_main} contains all $\Kt$-admissible lattices when $\Phi_\tau$ is of type $\mathsf{A}^{(2)}_{4}$.

\subsection*{Case: $\Phi_\tau$ is of type $\mathsf{A}_{2n-1}^{(2)}$ with $n\ge 3$}
\begin{lem} \label{lem:a2n12}
Suppose $\Phi_\tau$ is of type $\mathsf{A}_{2n-1}^{(2)}$ with $n\ge 3$ and let $\Xi$ be a $G(\sv_0)$-adapted sublattice of $\wl_{\sv_0}$ that does not satisfy \ref{al_doubles} at $\sv_0$. Then $\Xi$ is not $G(\sv_n)$-adapted, and therefore not $\Kt$-admissible. 
\end{lem}
\begin{proof}
For this affine root system $\Phi_\tau$, $G(\sv_0)$ is of type $\mathsf{C}_n$, with $n\ge 3$. As $\Xi$ does not satisfy \ref{al_doubles} at $\sv_0$, it follows that it satisfies \ref{al_Cn}, that is, $\Xi = \wl_{\sv_0}$. This implies that
\begin{equation}
\Xi \not\subset 2\wl_{\sv_n}. \label{eq:A2n12_lat}
\end{equation}  Indeed, $\w_1 \in \Xi$, $\aa_1 \in \oS(\sv_n)$ and $\la\aa_1^\vee, \w_1\ra = 1$, which is odd.  

If $n\ge 4$, then $G(\sv_n)$ is of type $\mathsf{D_n}$, which means, by \cref{fullrank}, that the only $G(\sv_n)$-adapted lattices are those satisfying \ref{al_doubles}. By \eqref{eq:A2n12_lat}, it follows that $\Xi$ is not $G(\sv_n)$-adapted.

If $n=3$, then $G(\sv_n) = G(\sv_3)$ is of type $\mathsf{A}_3$ and its Dynkin diagram is
\[ \dynkin[labels={0,2,1}, label macro/.code={\aa_{\drlap#1}}, edge length=.7cm] A{3} \]
By \eqref{eq:A2n12_lat}, $\Xi$ does not satisfy \ref{al_doubles} at $\sv_3$.  We show that it also doesn't satisfy \ref{al_An_odd} at $\sv_3$, which then implies by \cref{fullrank} that $\Xi$ is not $G(\sv_3)$-adapted, as there are no other adapted lattices for a group of type $\mathsf{A}_3$. Recall that $\w_1,\w_2,\w_3$ are the fundamental weights of $G(\sv_0)$. 
For the root system $\Phi_\tau$ of type $\mathsf{A}_{5}^{(2)}$, \cref{eq:corootequation} becomes
\[\aav_0 + \aav_1 + 2\aav_2 + 2\aav_3 = 0 \quad \text{or equivalently}\quad \aav_0 = -\aav_1-2\aav_2-2\aav_3.\]
Using this formula, one computes the matrix
\[
\begin{pmatrix}
\la \aa_0^\vee, \w_1\ra & \la \aa_1^\vee, \w_1\ra \\
\la \aa_0^\vee, \w_2\ra & \la \aa_1^\vee, \w_2\ra \\
\la \aa_0^\vee, \w_3\ra & \la \aa_1^\vee, w_3\ra \\
\end{pmatrix}
= 
\begin{pmatrix}
-1 & 1 \\
-2 & 0 \\
-2 & 0
\end{pmatrix}. 
\]
As the greatest common divisor of the $(2\times 2)$-minors of this matrix is $2$, the elementary divisors theorem tells us, that the coroots $\aa_0^\vee$ and $\aa_1^\vee$ are not part of a basis of the dual lattice $\Xi^*$. Consequently, $\Xi$ does not satisfy \ref{al_An_odd} at $\sv_3$.
\end{proof}

\Cref{lem:a2n12} and \cref{lem:aldoubles_adap} establish the following
\begin{lem}
Suppose $\Phi_\tau$ is of type $\mathsf{A}_{2n-1}^{(2)}$ with $n\ge 3$. The subgroups of $\fa$ that are $\Kt$-admissible are the lattices $\Xi$ with $2\oS(\sv_0) \subset \Xi \subset 2\lat_{\tau}$. 
\end{lem}

This shows that \cref{thm_main} contains all $\Kt$-admissible subgroups of $\fa$ when $\Phi_{\tau}$ is of type $\mathsf{A}_{2n-1}^{(2)}$ with $n\ge 3$.

\subsection*{Case: $\Phi_\tau$ is of type $\mathsf{D}_{n+1}^{(2)}$ with $n\ge 2$}

Here  $G(\sv_0)$ is of type $\mathsf{B}_n$. By \cref{fullrank} (and \cref{rem_after_fullrank}\ref{rem_after_fullrank_B2C2}), the $G(\sv_0)$-adapted lattices which do not satisfy \ref{al_doubles} at $\sv_0$ are 
\begin{align}
& \Z(\oS(\sv_0)^+\cup\{2\aa_n\}) \text{ and }\la \w_1,\w_2, \ldots, \w_{n-1},2\w_n\ra_\Z \text{ for all $n\ge 2$, and} \label{eq:dn122_gen} \\
&\text{in addition }\wl_{\sv_0} = \la \w_1,\w_2\ra_{\Z} \text{ when $n=2$}. \label{eq:dn122_2}
\end{align}
We will use that \cref{eq:delta_grad} and \cref{eq:corootequation} become 
\begin{align} 
&\aa_0+\aa_1+\cdots+\aa_n = 0,\quad \text{and}  \label{eq:delta_grad_dn122}\\
&\aav_0 +2\aav_1 + 2\aav_2 + \cdots +2\aav_{n-1} + \aav_n = 0 \label{eq:coroots_dn122}
\end{align}
for the affine root system $\Phi_{\tau}$ of type $\mathsf{D}_{n+1}^{(2)}$. 

\begin{lem} Suppose $\Phi_{\tau}$ is of type $\mathsf{D}_{n+1}^{(2)}$ with $n\ge 2$.
If $n$ is even, then $\Z(\oS(\sv_0)^+\cup\{2\aa_n\})$ is not $G(\sv_n)$-adapted and therefore not $\Kt$-admissible. 
\end{lem}
\begin{proof}
Set $\Xi = \Z(\oS(\sv_0)^+\cup\{2\aa_n\})$. Like $G(\sv_0)$, the group  $G(\sv_n)$ has Dynkin type $\mathsf{B}_n$. First note that $\la \aav_1,\aa_1+\aa_2\ra = 1$ (when $n=2$  too). As $\aa_1+\aa_2 \in \Xi$ and $\aa_1 \in \oS(\sv_n)$, this implies that $\Xi \not\subset 2\wl_{\sv_n}$ and so $\Xi$ does not satisfy \ref{al_doubles} at $\sv_n$.
Next we show that 
\begin{equation}
\label{eq:a01notinxi}
\aa_0 + \aa_1 \notin \Xi,
\end{equation}
 since this implies that $\Xi$ does not satisfy \ref{al_Bn_1} or \ref{al_Bn_2} at $\sv_n$ and that it does not satisfy \ref{al_Cn} when $n=2$. 
Since $\oS(\sv_0)^+\cup\{2\aa_n\}$ is a basis of $\Xi$, we know that $\aa_n \notin \Xi$. Because $n$ is even and because $\oS(\sv_0)^+ \in \Xi$, this implies \eqref{eq:a01notinxi} thanks to \cref{eq:delta_grad_dn122}. 
As we have shown that $\Xi$ cannot be any of the $G(\sv_n)$-adapted lattices listed in \cref{fullrank}, the lemma follows. 
\end{proof}

We now show that the lattice in \eqref{eq:dn122_2} is not $\Kt$-admissible.
\begin{lem}
Suppose $\Phi_\tau$ is of type $\mathsf{D}_{3}^{(2)}$. Then $\lat_\tau=\la \w_1,\w_2\ra_\Z$ is not $G(\sv_1)$-adapted and thefore not $\Kt$-admissible.  
\end{lem}
\begin{proof}
Here $G(\sv_1)$ is of type $\mathsf{A}_1 \times \mathsf{A}_1$. Writing $\wt_0$ and $\wt_2$ for the fundamental weights of $G(\sv_1)$ we find, using \cref{eq:coroots_dn122}, that $\Xi = \la -2\wt_0, \wt_2 -\wt_0\ra_\Z$, which is exacly the lattice we encountered in \cref{eq:ansl2sl2lat} in the proof of \cref{lem:c21}. We showed there that it follows from  \cref{pvsalg} that $\Xi$ is not $G(\sv_1)$-adapted. 
\end{proof}

Next we show that the remaining lattices in \cref{eq:dn122_gen} are $\Kt$-admissible. 

\begin{lem}
Suppose $\Phi_{\tau}$ is of type $\mathsf{D}_{n+1}^{(2)}$ with $n\ge 2$.
If $n$ is odd, then $\Z(\oS(\sv_0)^+\cup\{2\aa_n\})$ is a  $\Kt$-admissible subgroup of $\fa$. 
\end{lem}
\begin{proof}
Set $\Xi = \Z(\oS(\sv_0)^+\cup\{2\aa_n\})$ and . First we observe that \cref{eq:delta_grad_dn122} and the fact the $n$ is odd imply
\begin{align*}
2\aa_0 &= -2(\aa_1+\aa_2)-2(\aa_3+\aa_4)-\dots-2(\aa_{n-2}+\aa_{n-1})-2\aa_n; \quad \text{and}\\
\aa_0 + \aa_1 &= -(\aa_2+\aa_3)-(\aa_4+\aa_5)-\cdots-(\aa_{n-1}+\aa_{n}). 
\end{align*}
Consequently 
\begin{equation} \label{eq:2aaa01inxi}
2\aa_0, \aa_0+\aa_1 \in \Xi
\end{equation}
and $\Xi =  \Z(\oS(\sv_n)^+\cup\{2\aa_0\})$. 
This shows that $\Xi$ satisfies \ref{al_Bn_1} at $\sv_n$ and consequently is $G(\sv_n)$-adapted. 

We now fix $\ell \in \{1,2,\ldots,n-1\}$ and check that $\Xi$ is $G(\sv_\ell)$-adapted using \cref{pvsalg}. Observe that $G(\sv_\ell)$ has Dynkin type $\mathsf{B}_{\ell} \times \mathsf{B}_{n-l}$. To begin, we note that
\begin{equation} \label{eq:2aaa01notdbl}
\la \aav_k, \aa_{k-1}+ \aa_{k}\ra = \la \aav_k, \aa_{k}+ \aa_{k+1}\ra = 1\ \text{ for all $k \in \{1,2,\ldots,n-1\}$.}
\end{equation}
Since $\oS(\sv_0)^+\cup\{2\aa_n\} \subset \Xi$ it follows from  \cref{eq:2aaa01inxi,eq:2aaa01notdbl} that 
\[\SN(\Xi) = (\{2\aa_0,\aa_0+\aa_1, 2\aa_n\} \cup \oS(\sv_0)^+) \setminus \{\aa_{\ell-1}+ \aa_\ell, \aa_\ell + \aa_{\ell+1}\}.\]
An elementary, if somewhat lengthy computation then shows that $\oS_\Xi = \emptyset$ . Consequently the three conditions in \cref{pvsalg} are trivially justified. By \cref{lem:Gadap}\ref{lem:Gadap_item3} we can conclude that $\Xi$ is $\Kt$-admissible.
\end{proof}

\begin{lem}
Suppose $\Phi_{\tau}$ is of type $\mathsf{D}_{n+1}^{(2)}$ with $n\ge 2$.
Then $\Z\oS(\sv_0) = \la \aa_1, \aa_2, \ldots, \aa_{n}\ra_\Z$ is a $\Kt$-admissible subgroup of $\fa$. 
\end{lem}
\begin{proof}
By \cref{lem:Gadap}\ref{lem:Gadap_item3} it suffices to show that $\Z\oS(\sv_0)$ is $G(\sv_\ell)$-adapted for all $\ell \in \{1,2,\ldots, n\}$. We begin with $\ell = n$. Using \cref{eq:delta_grad_dn122} one directly sees, that 
\[\Z\oS(\sv_0) = \la \aa_0, \aa_1, \ldots, \aa_{n-1}\ra_{\Z} = \Z\oS(\sv_n).\]
Consequently, $\Z\oS(\sv_0)$ satisfies \ref{al_Bn_2} at $\sv_n$ and is therefore $G(\sv_n)$-adapted. We now fix $\ell \in \{1,2,\ldots,n-1\}$. Then $\overline{\Phi_\tau(\sv_\ell)}$ is of type $\mathsf{B}_\ell \times \mathsf{B}_{n-\ell}$ and $G(\sv_\ell) = G_1 \times G_2$, where $G_1$ has type $\mathsf{B}_\ell$ and set of simple roots $\{\aa_0, \aa_1, \ldots, \aa_{\ell-1}\}$, and $G_2$ has type $\mathsf{B}_{n-\ell}$ and set of simple roots $\{\aa_{\ell+1}, \aa_{\ell+2}, \ldots,\aa_n\}$. 
Using \cref{eq:delta_grad_dn122} once again, one directly sees that 
\(\Z\oS(\sv_0) =\Chi_1 \oplus \Chi_2\), where \[\Chi_1 = \la \aa_0, \aa_1, \ldots, \aa_{\ell-1}\ra_\Z \quad\text{and}\quad \Chi_2=\la \aa_{\ell+1}, \aa_{\ell+2}, \ldots, \aa_{n}\ra_\Z .\]
Consequently, $\Chi_1$ is $G_1$-adapted and $\Chi_2$ is $G_2$-adapted by \ref{al_Bn_2} and therefore $\Chi$ is $G(\sv_\ell)$-adapted. 
\end{proof}

We have shown the following
\begin{lem}
If  $\Phi_\tau$ is of type $\mathsf{D}_{n+1}^{(2)}$ with $n\ge 2$, then the $\Kt$-admissible subgroups of $\fa$ are:  $2\Z\oS(\sv_0)$, $2\lat_{\tau}$ and $\Z\oS(\sv_0)$ and, in addition,  $\Z(\oS(\sv_0)^+\cup\{2\aa_n\})$ when $n$ is odd.
\end{lem}
This justifies entries (\ref{spineven}) and (\ref{spinodd}) in \cref{thm_main} and shows that \cref{thm_main} contains all $\Kt$-admissible subgroups of $\fa$ when $\Phi_\tau$ is of type $\mathsf{D}_{n+1}^{(2)}$ with $n\ge 2$. 

\def\cprime{$'$} \def\cprime{$'$} \def\cprime{$'$} \def\cprime{$'$}
  \def\cprime{$'$}
\providecommand{\bysame}{\leavevmode\hbox to3em{\hrulefill}\thinspace}
\providecommand{\MR}{\relax\ifhmode\unskip\space\fi MR }
% \MRhref is called by the amsart/book/proc definition of \MR.
\providecommand{\MRhref}[2]{%
  \href{http://www.ams.org/mathscinet-getitem?mr=#1}{#2}
}
\providecommand{\href}[2]{#2}

\end{document}